\magnification 1200

\centerline{\bf Trees of metric compacta and trees of manifolds}
\centerline{Jacek \'Swi\c atkowski
\footnote{*}{This work was supported by the Polish Ministry 
of Science and Higher Education (MNiSW), grant N201 541738,
and by the Polish National Science Centre (NCN), grant 2012/06/A/ST1/00259.}}


\bigskip
\bigskip
\centerline{\bf Introduction}

\bigskip
We present a construction, called the
{\it limit of a tree system of spaces} (or, less formally,
a {\it tree of spaces}). The construction is designed to
produce compact metric spaces that resemble fractals, 
out of more regular spaces,
such as closed manifolds, compact polyhedra, compact Menger
manifolds, etc. Such spaces are potential candidates to be
homeomorphic to ideal boundaries of infinite groups.

A very special case of this construction, {\it trees of manifolds}
(known also as {\it Jakobsche spaces}), 
has been studied in the literature
(see [AS], [J1], [J2], [St]).
We present here a different approach, much more general, and, as we believe,
much more convenient for establishing various basic properties 
of the resulting spaces, in a more general setting.
Already in the case of trees of manifolds, using this approach 
we clarify, correct and extend so far known results and properties.

\medskip\noindent
{\it A wider context for the results of the present paper.}

Our motivation for dealing with the general construction as presented in
this paper comes from an attempt to understand which topological spaces
are boundaries of hyperbolic groups. This problem, stated e.g. in [KK]
as Question A,  
remains widely open.
An overview of the limited knowledge concerning this problem
can be found in Section 17 of [KB] or in the introduction in [P\'S].

The present paper initiates a larger research project,
investigated by the author, concerning the above problem. 
We briefly outline the aims and expected
lines of further development in this project.

\item{$\bullet$}
In a paper under preparation we describe a vast class of topological spaces
called {\it trees of polyhedra}. These spaces are compact, metrizable,
have finite topological dimension, and typically they are "wild"
(e.g. they are usually not locally contractible, and hence not ANR). 
They are obtained as limits of some
tree systems, and depend uniquely up to homeomorphism
on certain finite data, part of which is a finite collection of compact
polyhedra. Thus, the spaces are given not by a universal characterization
in terms of a list of properties, but rather by a sort of ``presentation''
(similar in spirit to a presentation of a group in combinatorial group
theory). Typically, the same space has many distinct ``presentations'',
and clarification of the relationship between such ``presentations''
will be one of the challenges.

\item{$\bullet$}
Next part of the project consists of identifying boundaries
of various classes of groups as explicit trees of polyhedra.
We have formulated several conjectures in this direction,
based on some results from the literature (see Remark 2.C.8), 
and on our new
partial results. One of these conjectures deals with Gromov boundaries
of all groups obtained by any procedure of strict hyperbolization 
(e.g. the one described in [CD]).
Another conjecture concerns CAT(0) boundaries of a large class
of right angled Coxeter groups, namely the ones that enjoy some
partial hyperbolicity property, weaker than word hyperbolicity.
We have already confirmed these conjectures in some cases
when the boundaries are trees of manifolds (in arbitrary dimension),
trees of graphs, or trees of various 2-dimensional polyhedra
which can be recognized as the Menger curve or the Sierpi\'nski curve.
It seems that boundaries of numerous (maybe even most of) other
groups so far studied in the literature are trees of polyhedra.

\item{$\bullet$}
A question arises, for which hyperbolic groups their Gromov boundaries
are not trees of polyhedra. 
The following examples seem to belong to this class:

\itemitem{(1)} groups whose Gromov boundaries are Menger compacta 
of dimension $\ge2$
(see [DO] for examples of such groups, with boundaries of dimension
2 and 3); we suspect that Menger compacta 
satisfy stronger disjoint
disk properties than any trees of polyhedra of the same 
topological dimension;

\itemitem{(2)} 7-systolic groups, as defined in [JS], with boundaries 
of dimension $\ge3$; by [Sw1], Gromov boundaries of 7-systolic groups
contain no copy of the 2-disk, which cannot happen for a tree
of polyhedra in dimensions $\ge3$;

\itemitem{(3)} topologically rigid hyperbolic groups, examples of which
have been constructed in [KK]; we expect that trees of polyhedra
always admit homeomorphisms whose dynamics is different from
that occuring for the induced action of a group on its boundary.

\item{} It seems to be known that boundaries of all hyperbolic groups
fall in the class of spaces called Markov compacta
(as defined e.g. in [Dr]). We work on showing that trees of polyhedra
coincide with a subclass of Markov compacta of some rather simple form.
This suggests a possibility to introduce certain notion of
degree of complexity for Markov compacta, with the lowest degree
corresponding to trees of polyhedra. A big challenge is to explore
more fully the territory of hyperbolic groups with higher degree
boundaries (i.e. boundaries which are not trees of polyhedra).
For example, one can ask for new strict hyperbolization procedures,
resulting with groups of higher degree boundaries. One can also ask
for explicit description of such boundaries, perhaps starting
with degree just above the lowest one.

\medskip\noindent
{\it The content of the paper.}

In Part 1 we introduce the notions of a tree system of metric compacta
(Section 1.B) and its limit (Section 1.C). We then prove that such limit
is always a compact metrizable space (Section 1.D)
by showing that it is homeomorphic to the limit of some inverse
system naturally associated to a tree system. In Section 1.E
we introduce the notion of isomorphism of tree systems.
We also present a class of natural examples - dense tree systems
of closed (topological)
manifolds, and a class of spaces obtained as their limits - 
the trees of manifolds $M$. As we show in Part 2 of the paper,
trees of manifolds $M$ coincide with the spaces studied earlier
by W. Jakobsche in [J2] (for $M$ oriented) and by P. Stallings in [St]
(for non-orientable $M$). Our exposition of the case of non-orientable 
manifolds $M$, in Subsection 1.E.4, concerns all topological manifolds
(and not only PL ones, as in [St]), and it frees the description from
certain unconvenient and unnecessary condition present in Stallings'
approach (see Remark 1.E.4.3). 

In Part 2 of the paper we 
exhibit some other useful inverse systems associated to tree systems.
In Section 2.A we introduce some additional data, called {\it extended
system of spaces and maps}, necessary to produce such inverse systems.
In Section 2.B we show that inverse limits of these new inverse systems
canonically coincide with limits of the corresponding tree systems.
Section 2.C deals with a subclass of tree systems called peripherally ANR,
and describes inverse systems of particularly nice form associated to 
such tree systems. This allows us to relate our construction of limit
of a tree system, in the case of dense tree systems of manifolds,
to some earlier constructions from the literature, notably the construction
of Jakobsche [J2].
In Section 2.D we apply associated inverse systems to provide estimates
from above for the topological dimension of limits of tree systems.
We indicate some general cases in which these estimates are sharp.

In Part 3 of the paper we introduce some natural and useful operations on
tree systems. In Section 3.A we describe an operation of consolidation,
by which  
the spaces appearing in the initial system are merged together
into bigger spaces, constituting naturally a new tree system.
We show that this operation does not affect the limit.
As an application, 
we derive equalities (up to homeomorphism)
between the limits of various different tree systems of manifolds.
We also include a correction to the main result of H. Fischer in [Fi],
in which he identifies boundaries of certain right angled Coxeter groups
as some explicit trees of manifolds $M$ (see Theorem 3.A.3 in the text,
and a comment after its statement).
In Section 3.B we show how to decompose a compact metric space $X$
into pieces which form a tree system, so that $X$ naturally
coincides with the limit of this tree system. We also show how to use
such decompositions to determine limits of certain tree systems.
In Section 3.C we introduce subdivision of a tree system,
an operation opposite to consolidation, which generalizes the
operation of decomposition from Section 3.B.
In Section 3.D we apply operations of subdivision and consolidation
to study orbits under homeomorphisms in trees 
of manifolds. The presented method is potentially applicable
to more general tree systems. 
Finally, in Section 3.E we apply operations of consolidation and subdivision
to provide a much more flexible description of certain trees of manifolds
than those so far present in the literature
(see Theorem 3.E.2 and Corollary 3.E.4 in the text). 
In the companion paper [Sw2],
we use this description to identify ideal boundaries of many groups
as trees of manifolds. In particular, we show in [Sw2] that trees of manifolds
in arbitrary dimension appear as Gromov boundaries of certain
hyperbolic groups.

\medskip\noindent
{\it Acknowledgments.}

The author thanks the referee for many suggestions leading to the
improvement of the present paper, and especially for inventing and communicating to the author the idea
of the standard inverse system associated to a tree system 
(presented in Section 1.D), and its 
usefulness in the whole of the paper.

\vfill
\break

\bigskip
\centerline{\bf 1. Tree systems and their limits.}

\bigskip\noindent
{\bf 1.A Some terminology and notation concerning trees.}

\medskip
Trees under our consideration will be
usually countable infinite, and locally infinite. We donote by $V_T$ the set
of all verticess of a tree $T$, and by $O_T$ the set of its all oriented
edges. For any $e\in O_T$, we denote by $\alpha(e),\omega(e)$ the initial and
the terminal vertex of $e$, respectively. We also denote by $\bar e$
the same geometric edge as $e$, but oppositely oriented.
For any $t\in V_T$, we denote by $N_t=\{e\in O_T:\alpha(e)=t  \}$ the set of 
all oriented edges of $T$ with initial vertex $t$.

We denote the combinatorial (embedded) paths in $T$ as sequences
of consecutive vertices $[t_0,t_1,\dots,t_m]$, or as sequences
$[e_1,\dots,e_m]$ of consecutive oriented edges, or shortly by
$[t,s]$ if $t$ and $s$ are the ends of the path.
An infinite combinatorial path $[t_0,t_1,\dots]$ in $T$ is called a $ray$,
and denoted usually by $\rho$. We denote by $\rho(0)$ the initial
vertex $t_0$, and by $e_1(\rho)$ the initial oriented edge $(t_0,t_1)$
of a ray $\rho$.

We denote by $E_T$ the set of {\it ends} of $T$, i.e. the set of
equivelence classes of rays in $T$ with respect to the relation
of coincidence except possibly at some finite initial parts.
We denote the end determined by a ray $\rho$ as $[\rho]$.

Let $S$ be a subtree of $T$. We then distinguish the set
$$
N_S=\{ e\in O_T:\alpha(e)\in V_S\,\hbox{ and }\,\omega(e)\notin V_S  \}.
$$
Note that, in case when $S$ is reduced to a single vertex $t$,
this notation agrees with the earlier introduced notation 
for the set $N_t$.
A finite subtree of $T$ will be usually denoted by $F$,
and we shall consider the poset $({\cal F}_T,\subset)$ of all
finite subtrees of $T$.

\bigskip\noindent
{\bf 1.B Tree systems of spaces.}

\medskip
Recall that a family of subsets of a compact metric space is {\it null}
if for each $\epsilon>0$ all but finitely many of these subsets have
diameters less than $\epsilon$. 

\medskip\noindent
{\bf 1.B.1 Definition.}
A {\it tree system of metric compacta} is a tuple 
$\Theta=(T,\{ K_t \},\{\Sigma_e\},\{\phi_e\})$ such that:
\itemitem{(TS1)} $T$ is a countable tree;
\itemitem{(TS2)} to each $t\in V_T$ there is associated a compact metric 
space $K_t$;
\itemitem{(TS3)} to each $e\in O_T$, there is associated 
a nonempty compact
subset $\Sigma_e\subset K_{\alpha(e)}$, and a homeomorphism 
$\phi_e:\Sigma_{e}\to\Sigma_{\bar e}$ such that 
$\phi_{\bar e}=\phi_e^{-1}$;
\itemitem{(TS4)} for each $t\in V_T$ the family 
$\Sigma_e:e\in N_t$
is null and consists of pairwise disjoint sets.

\noindent
We call $T$ the {\it underlying tree}, $K_t:t\in V_T$ the {\it constituent
spaces}, $\Sigma_e:e\in O_T$ the {\it peripheral subspaces},
and $\phi_e:e\in O_T$ the {\it connecting maps} of the tree system
$\Theta$.

\medskip\noindent
{\bf Remark.} In future applications of tree systems of spaces
we will often additionally require 
that for any $t\in V$ the family $\Sigma_e:e\in N_t$
is {\it dense} in the space $K_t$ (which means that the union of this
family is a dense subset). However, to establish many basic properties
of tree systems of metric compacta we do not need this requirement.

\medskip\noindent
{\bf 1.B.2 Example: tree system of manifolds.}
Let $T$ be a countable tree. Let $M_t:t\in V_T$ be a family of
closed manifolds of the same dimension $n$.
For each $t\in V_T$ and each $e\in N_t$ let $\Delta_e$ 
be a collared $n$-dimensional disk
embedded in $M_t$, and suppose that each 
${\cal D}_t=\{ \Delta_e:e\in N_t \}$ is a null family of pairwise disjoint
subsets of $M_t$.

For each $t\in V_T$ put
$$
K_t=M_t\setminus
\bigcup\{ \hbox{int}(\Delta_e):e\in N_t \}.
$$
This defines a family $\{K_t\}$ in a tree system of manifolds.

For each $e\in O_T$, put $\Sigma_e=\partial\Delta_e$, 
and note that $\Sigma_e\subset K_{\alpha(e)}$.
For each $e\in O_T$ consider also a homeomorphism 
$\phi_e:\Sigma_{\alpha(e)}\to \Sigma_{\omega(e)}$
between the corresponding $(n-1)$-spheres, so that
$\phi_{\bar e}=\phi_e^{-1}$ for each $e$. 
A {\it tree system of manifolds} is a tuple 
${\cal M}=(T,\{ K_t \}, \{ \Sigma_e \}, \{ \phi_e \})$ as described
above.

Intuitively, such a system $\cal M$ may be viewed as a
pattern for a connected sum operation applied at the same
time to a countable (in general infinite) family of 
closed manifolds.

\bigskip\noindent
{\bf 1.C Limit of a tree system of spaces.}

\medskip
We now describe the {\it limit} of a tree system 
$\Theta=(T,\{ K_t \},\{\Sigma_e\},\{\phi_e\})$, denoted $\lim\Theta$,
starting with its description as a set. Denote by $\#\Theta$
the quotient 
$$(\bigsqcup_{t\in V_T}K_t)/\sim$$
of the disjoint union of the sets $K_t$ by the equivalence relation $\sim$
induced by the equivalences $x\sim\phi_e(x)$ for any $e\in O_T$ and any
$x\in\Sigma_e$. This set may be viewed as obtained from the family
$K_t:t\in V$ as a result of gluings provided by the maps $\phi_e$. 
Observe that any set $K_t$ canonically injects
in $\#\Theta$.
Define $\lim\Theta$ to be the disjoint union $\#\Theta\cup E_T$,
where $E_T$ is the set of ends of $T$.

To put appropriate topology on the set $\lim\Theta$ we need some
terminology. Given a family $\cal A$ of subsets in a set $X$, we say that
$U\subset X$ is {\it saturated} with respect to $\cal A$ 
(shortly, $\cal A$-{\it saturated}) if for any $A\in{\cal A}$
we have either $A\subset U$ or $A\cap U=\emptyset$.

For any finite subtree $F$ of $T$, denote by $V_F,O_F$ the sets of vertices
and of oriented edges in $F$, respectively.
Denote also by 
$$
\Theta_F=(F,\{ K_t:t\in V_F \},
\{ \Sigma_e:e\in O_F \},
\{ \phi_e:e\in O_F \})
$$ 
the tree system of spaces $\Theta$ {\it restricted}
to $F$. The set $K_F:=\#\Theta_F$, equipped with the natural quotient topology
(with which it is clearly a metrizable compact space), will be called
the {\it partial union} of the system $\Theta$ related to $F$.
Observe that any partial union $K_F$ is canonically
a subset in $\#\Theta$, and thus also in $\lim\Theta$.

For any finite subtree $F$ put ${\cal A}_F=\{ \Sigma_e:e\in N_F \}$,
and view the elements of this family as subsets in the partial union 
$K_F$. Observe that the family ${\cal A}_F$ consists of pairwise disjoint
compact sets, and that it is null. Let $U\subset K_F$ be a subset which is
saturated with respect to ${\cal A}_F$.
Put $N_U:=\{ e\in N_F:\Sigma_e\subset U \}$ and 
$$
D_U:=\{ t\in V_T:[\omega(e),t]\cap V_F=\emptyset
\hbox{ for some }e\in N_U \}
$$
(i.e. $D_U$ is the set of these vertices $t\in V_T\setminus V_F$
for which the shortest path in $T$ connecting $F$ to $t$
starts with an edge $e\in N_U$).
Here we mean in particular that the vertices $\omega(e):e\in N_U$
all belong to $D_U$.

Again, for any finite subtree $F$ of $T$ let $R_F$ be the set of
all rays $\rho$ in $T$ with $\rho(0)\in V_F$, and with all other
vertices outside $V_F$. Note that the map $\rho\to[\rho]$
is then a bijection from $R_F$ to the set $E_T$ of all ends of $T$.
Given a subset $U\subset K_F$ as above (i.e. saturated with respect to
${\cal A}_F$), put $R_U:=\{ \rho\in R_F:e_1(\rho)\in N_U \}$
and $E_U:=\{ [\rho]\in E_T:\rho\in R_U \}$.

Finally, for any subset $U\subset K_F$ as above, put
$$
G(U):=U\cup(\bigcup_{t\in D_U}K_t)\cup E_U,
$$
where all summands in the above union are viewed as subsets of $\lim\Theta$;
clearly, $G(U)$ is then also a subset in $\lim\Theta$.
We consider the topology in the set $\lim\Theta$ given by the basis $\cal B$
consisting of all sets $G(U)$, for all finite subtrees $F$ in $T$,
and all open subsets $U\subset K_F$ saturated with respect to
${\cal A}_F$.
The family $\cal B$ satisfies the axioms
for a basis of topology because intersection of any two sets from $\cal B$
is again in $\cal B$, as it can be quickly deduced from the following two
easy observations:
\item{(1)} if $U,U'\subset K_F$ are open and ${\cal A}_F$-saturated 
then $U\cap U'$ is also ${\cal A}_F$-saturated, and 
$G(U)\cap G(U')=G(U\cap U')$;
\item{(2)} if $U\subset K_F$ is open and ${\cal A}_F$-saturated,
then for any finite $F'\supset F$ the set $U':=G(U)\cap K_{F'}$
is ${\cal A}_{F'}$-saturated, and $G(U')=G(U)$.

\medskip\noindent
{\bf 1.C.1 Proposition.} {\it For any tree system $\Theta$ of metric compacta
the limit $\lim\Theta$, with topology given by the above described basis $\cal B$, 
is a metrizable compact space. Moreover, for each
finite subtree $F$ of $T$
the canonical inclusion (as a set) of the partial union 
$K_F=\#\Theta_F$ in $\lim\Theta$
is a topological embedding.
In particular, for any $t\in V_T$ the canonical inclusion of $K_t$
in $\lim\Theta$ is a topological embedding.
Finally, both families of subsets $\{ K_t:t\in V_T \}$
and $\{\Sigma_e:e\in O_T  \}$ in $\lim\Theta$ are null
(with respect to any metric compatible with the topology).}

\medskip
We present a proof of Proposition 1.C.1 in the next section,
after showing Proposition 1.D.1, which we need in this proof.

\medskip\noindent
{\bf 1.C.2 Remark.} If we skip the assumption that each of the families
of  subsets $\{ \Sigma_e:e\in N_t \}$ in $K_t$, for any $t\in V_T$, is null,
then the above description still defines a compact topological space
$\lim\Theta$, but then the limit is in general not Hausdorff.
Moreover, canonical inclusions $K_t\to\lim\Theta$ are then in general
not embeddings (though they are always continuous).

\bigskip\noindent
{\bf 1.D The standard inverse system associated to a tree system.}

\medskip
We describe some inverse system of metric compacta naturally
associated to a tree system $\Theta$,
and show that the inverse limit of this system is cononically
homeomorphic to the limit $\lim\Theta$. We use this fact to prove
Proposition 1.C.1, and to derive few further properties of limits
of tree systems.

Given a tree system
$\Theta=(T,\{ K_t \},\{ \Sigma_e \},\{ \phi_e \})$, 
we proceed using the notation as in Section 1.C.
For any finite subtree $F\subset T$,
denote by $K_F^*=K_F/{\cal A}_F$ the qotient of $K_F$
in which all the subsets ${\cal A}_F$ are shrinked to points.
More precisely, let ${\cal D}_F$ be the decomposition of $K_F$
consisting of the sets in ${\cal A}_F$ and the singletons
from the complement $K_F\setminus\bigcup{\cal A}_F$,
and let $K_F^*=K_F/{\cal D}_F$. Since the family ${\cal A}_F$ is null,
the decomposition ${\cal D}_F$ is upper semicontinuous,
and hence the quotient $K_F^*$ is metrizable (see [Dav], Proposition 3
on p. 14 and Proposition 2 on p.13). Consequently, $K_F^*$ is a metric
compactum, and we call it the {\it reduced partial union of $\Theta$}
related to $F$. We denote by $q_F:K_F\to K_F^*$ the quotient map
resulting from the above description of $K_F^*$.

For any pair $F_1\subset F_2$ of finite subtrees of $T$, define
a map $f_{F_1F_2}:K_{F_2}^*\to K_{F_1}^*$ as follows.
For each edge $e\in N_{F_1}\cap O_{F_2}$ denote by ${\cal V}_e$
the set of all vertices $s\in V_{F_2}\setminus V_{F_1}$
such that the shortest path in $T$ connecting $F_1$ with $s$
starts with $e$. Then the family ${\cal V}_e:e\in N_{F_1}\cap O_{F_2}$
is a partition of $V_{F_2}\setminus V_{F_1}$, and each ${\cal V}_e$
is the vertex set of a subtree of $F_2$ which we denote $S_e$.
Viewing each $K_{S_e}$ canonically as a subset in $K_{F_2}$, note
that it is ${\cal A}_{F_2}$-saturated. Thus,
the corresponding subset $q_{F_2}(K_{S_e})$ in $K_{F_2}^*$
is well defined and closed. Moreover, by shrinking each of the subsets
$q_{F_2}(K_{S_e}):e\in N_{F_1}\cap O_{F_2}$ to a point we get
a quotient of $K_{F_2}^*$ which is canonically homeomorphic
to (and which we identify with) $K_{F_1}^*$. Take the
corresponding quotient map as $f_{F_1F_2}$, and observe
that this map is continuous and surjective.

Given any finite subtrees $F_1\subset F_2\subset F_3$ of $T$,
it is not hard to see that $f_{F_1F_2}\circ f_{F_2F_3}=f_{F_1F_3}$.
Consequently, the system
$$
{\cal S}_\Theta=(\{ K_F^*:F\in{\cal F}_T \},\{ f_{FF'}:F\subset F' \})
\leqno{(1.\hbox{D}.1)}
$$
is an inverse system of metric compacta over the poset ${\cal F}_T$
of all finite subtrees of $T$. We call it the {\it standard inverse 
system associated to} $\Theta$.

\medskip
We now turn to describing a natural map 
$\beta:\lim\Theta\to\lim_{\longleftarrow}{\cal S}_\Theta$
(initially as a map between the sets, but then we show that it is
a homeomorphism).
We split the description into two parts, accordingly with the partition
$\lim\Theta=\#\Theta\sqcup E_T$.

To describe $\beta$ on the subset $\#\Theta\subset\lim\Theta$,
for each finite subtree $F\subset T$ consider the map 
$f_{F\#}:\#\Theta\to K_F^*$ defined similarly to the maps $f_{FF'}$, 
as follows.
For each edge $e\in N_{F}$ denote by ${\cal V}_e$
the set of all vertices $s\in V_T\setminus V_{F}$
such that the shortest path in $T$ connecting $F$ with $s$
starts with $e$. Then the family ${\cal V}_e:e\in N_{F}$
is a partition of $V_T\setminus V_{F}$, and each ${\cal V}_e$
is the vertex set of a subtree of $T$ which we denote $S_e$.
For each $e\in N_F$,
denote by $\Theta_{S_e}$ the tree system obtained by restricting 
$\Theta$ to $S_e$, and
view $\#\Theta_{S_e}$ canonically as a subset in $\#\Theta$.
Note that $\#\Theta$ splits as
$$
\#\Theta=(K_F\setminus\bigcup_{e\in N_F}\Sigma_e)\sqcup
\bigsqcup_{e\in N_F}\#\Theta_{S_e}.
$$
Viewing $K_F$ as a subset in $\#\Theta$, for any 
$x\in K_F\setminus\bigcup_{e\in N_F}\Sigma_e$
we put $f_{F\#}(x):=x$. For any $e\in N_F$ and any  $x\in\#\Theta_{S_e}$
we put $f_{F\#}(x):=q_F(\Sigma_e)\in K_F^*$, where the latter is a single point
by definition of $K_F^*$.
A straightforward verification shows that for any finite subtrees
$F\subset F'$ in $T$ we have $f_{FF'}\circ f_{F'\#}=f_{F\#}$.
Thus the family $f_{F\#}:F\in{\cal F}_T$ induces
a well defined map $f_\#:\#\Theta\to\lim_{\longleftarrow}{\cal S}_\Theta$,
and we put $\beta(x):=f_\#(x)$ for any $x\in\#\Theta$.

To describe $\beta$ on the subset $E_T\subset\lim\Theta$,
consider any $x\in E_T$. For any finite subtree $F\subset T$ let $e^F_x$
be the first oriented edge in the "shortest" path in $T$ connecting $F$ to $x$
(i.e. in the unique ray in $T$ starting at a vertex of $F$,
passing through no other vertex of $F$, and 
representing $x$). Put $x(F)=q_F(\Sigma_{e^F_x})\in K_F^*$.
It is not hard to observe that the tuple $(x(F))_{F\in{\cal F}_T}$
is a thread of the inverse system ${\cal S}_\Theta$, i.e. an element
of the inverse limit $\lim_{\longleftarrow}{\cal S}_\Theta$,
which we denote $\xi_x$. We then put $\beta(x):=\xi_x$.

\medskip\noindent
{\bf 1.D.1 Proposition.}
{\it The above described map 
$\beta:\lim\Theta\to\lim_{\longleftarrow}{\cal S}_\Theta$
is a homeomorphism.}

\medskip\noindent
{\bf Proof:}
We will first show that $\beta$ is a bijection, and then that both $\beta$
and $\beta^{-1}$ are continuous.

\medskip\noindent
{\it $\beta$ is injective.}

Consider first any two distinct points $x,y\in\#\Theta$.
Let $F$ be any finite subtree of $T$ such that
$x,y\in K_F\setminus\bigcup_{e\in N_F}\Sigma_e$
(which obviously exists). By definition of the map $f_{F\#}$, we get 
$f_{F\#}(x)\ne f_{F\#}(y)$, and this implies $f_\#(x)\ne f_\#(y)$,
and hence also $\beta(x)\ne\beta(y)$.

Second, consider any two distinct $x,y\in E_T$.
Obviously, there is a finite subtree $F\subset T$ such that the 
oriented edges, $e_x^F$ and $e_y^F$ are distinct.
Since then $x(F)\ne y(F)$, we get $\xi_x\ne\xi_y$,
and thus $\beta(x)\ne\beta(y)$.

To finish the proof of injectivity, we need to show that for any 
$x\in\#\Theta$ and any $y\in E_T$, we have $\beta(x)\ne\beta(y)$.
To do this, for any $F\in{\cal F}_T$ denote by $P_F$ the set of all
points in $K_F^*$ of form $q_F(\Sigma_e):e\in N_F$.
Observe that $\beta(y)=(y(F))_{F\in{\cal F}_T}$ has the property
that for each $F\in{\cal F}_T$ we have $y(F)\in P_F$.
On the other hand, denoting $\beta(x)=(x_F)_{F\in{\cal F}_T}$,
we have $x_F=f_{F\#}(x)$ for any $F\in{\cal F}_T$.
Taking $F_0$ such that 
$x\in K_{F_0}\setminus\bigcup_{e\in N_{F_0}}\Sigma_e$,
we get $x_{F_0}\notin P_{F_0}$, and thus $\beta(x)\ne\beta(y)$,
as required.

\medskip\noindent
{\it $\beta$ is surjective.}

Let $x_0=(x_F)_{F\in{\cal F}_T}$ be an arbitrary point of 
$\lim_{\longleftarrow}{\cal S}_\Theta$. Suppose first that for some
$F_0\in{\cal F}_T$ we have $x_{F_0}\notin P_F$.
We may then view $x_{F_0}$ as a point of 
$K_{F_0}\setminus\bigcup_{e\in N_{F_0}}\Sigma_e$.
Observe that for any finite $F\supset F_0$ there is unique $z_F\in K_F^*$
such that $f_{F_0F}(z_F)=x_{F_0}$, and in fact we have
$z_F=q_F(x_{F_0})=f_{F\#}(x_{F_0})$, under canonical inclusions
$x_{F_0}\in K_{F_0}\subset K_F\subset\#\Theta$. From this it easily follows
that, viewing $x_{F_0}$ as an element of $\#\Theta$, we have
$x_0=\beta(x_{F_0})$.

In the remaining case, for each $F\in{\cal F}_T$ we have
$x_F\in P_F$. As we will see, in this case the tree $T$ is necessarily
unbounded. For any $F\in{\cal F}_T$ denote by $e_F$ the edge
in $N_F$ corresponding to $x_F\in P_F$ (i.e. such that
$x_F=q_F(\Sigma_{e_F})$), and put $t_F:=\alpha(e_F)\in V_F$.
Consider any increasing sequence $\sigma=(F_i)_{i\ge1}$ of finite subtrees
of $T$ such that $\bigcup_{i\ge1}F_i=T$ (we will call any such sequence
an {\it exhausting sequence}). It is not hard to realize that the 
associated sequence $(t_{F_i})_{i\ge1}$, after deleting potential repetitions
of subsequent terms, is necessarily infinite, and 
determines a ray in $T$ starting at $t_{F_1}$
and passing through all vertices $t_{F_i}$ as well as through
all oriented edges $e_{F_i}$ 
(and possibly through some other vertices and edges). We denote this ray
by $\rho_\sigma$.
Given any two exhausting sequences $\sigma=(F_i)_{i\ge1}$ and
$\sigma'=(F_i')_{i\ge1}$ of finite subtrees in $T$, we obviously 
have the following property: 
{\it for any $i\ge1$ there is $j\ge1$ such that $F_i\subset F_j'$
and $F_i'\subset F_j$.} From this property one deduces that the rays
$\rho_\sigma$ and $\rho_{\sigma'}$ eventually coincide.
It follows that in the considered case the point
$x_0=(x_F)_{F\in{\cal F}_T}$ uniquely determines an end $\xi_0\in E_T$.
It is also clear that $\beta(\xi_0)=x_0$,
which completes the proof.

\break

\medskip\noindent
{\it $\beta$ and $\beta^{-1}$ are continuous.}

Recall that, by definition, the inverse limit 
$\lim_{\longleftarrow}{\cal S}_\Theta$ is a subspace in the topological
product $\prod_{F\in{\cal F}_T}K_F^*$. Given any open subset
$U\subset K_{F_0}^*$, for any $F_0\in{\cal F}_T$, put
$$
G_{\cal S}(U):=(U\times\prod_{F\in{\cal F}_T\setminus\{F_0\}}K_F^*)\,
\cap \, \lim_{\longleftarrow}{\cal S}_\Theta.
$$
Note that $G_{\cal S}(U)$ is an open subset of 
$\lim_{\longleftarrow}{\cal S}_\Theta$, and that all subsets of this form
constitute a subbasis of the topology of 
$\lim_{\longleftarrow}{\cal S}_\Theta$.

Put $U':=q_{F_0}^{-1}(U)$ and note that, as $U$ runs through all
open subsets of $K_{F_0}^*$, $U'$ runs through all open 
${\cal A}_{F_0}$-saturated subsets of $K_{F_0}$.
Moreover, a direct observation shows that 
$$
\beta^{-1}(G_{\cal S}(U))=G(U'), \leqno{(1.\hbox{D}.2)}
$$
where $G(U')$ is an element of the basis $\cal B$ (of the topology
in $\lim\Theta$) described
in Section 1.C.
Since, by what was said above, (1.D.2) implies that both $\beta$
and $\beta^{-1}$ are continuous,
this completes the proof of Proposition 1.D.1.

\medskip\noindent
{\bf Proof of Proposition 1.C.1:}
We apply Proposition 1.D.1. Since $\lim\Theta$ is homeomorphic
to the inverse limit of an inverse system of metric compacta, it is itself
a compact metrizable space, which yields the first assertion.

To see the second assertion, for each $F_0\in{\cal F}_T$ consider
the subposet ${\cal F}_T^{F_0}$ consisting of all subtrees $F\in{\cal F}_T$
that contain $F_0$, and note that this subposet is cofinal with ${\cal F}_T$.
Denote by ${\cal S}_\Theta^{F_0}$ the restriction
of the inverse system ${\cal S}_\Theta$ to this subposet,
and note that we have a canonical identifications of the limits 
$\lim_{\longleftarrow}{\cal S}_\Theta^{F_0}=
\lim_{\longleftarrow}{\cal S}_\Theta$. For each $F\in{\cal F}_T^{F_0}$
consider the map $h_F^{F_0}:K_{F_0}\to K_F^*$ given as the composition
of the inclusion map $K_{F_0}\to K_F$ and the map $q_F:K_F\to K_F^*$.
Note that the family of maps $h_F^{F_0}:F\in{\cal F}_T^{F_0}$
gives a morhism $K_{F_0}\to{\cal S}_\Theta^{F_0}$,
and thus induces a continuous map 
$h^{F_0}:K_{F_0}\to\lim_{\longleftarrow}{\cal S}_\Theta$.
Since it is not hard to observe that, under identification of
$\lim_{\longleftarrow}{\cal S}_\Theta$ with $\lim\Theta$ through $\beta$,
the map $h^{F_0}$ coincides with the inclusion map 
$K_{F_0}\to\lim\Theta$, it follows that the latter map is an embedding.

To see the last assertion,
consider an auxilliary increasing sequence $(F_i)_{i\ge1}$ 
of finite subtrees of $T$
such that $T=\bigcup_{i\ge1} F_i$. Obviously, this sequence yields
a cofinal subposet of the poset ${\cal F}_T$. Denote by ${\cal S}'$
the inverse sequence obtained by restricting ${\cal S}_\Theta$
to this subposet, and note that there is a canonical identification
$\lim_{\longleftarrow}{\cal S}'=\lim_{\longleftarrow}{\cal S}_\Theta$.
The map $\beta$ is then naturally viewed as a map 
$\lim\Theta\to\lim_{\longleftarrow}{\cal S}'$.
Call a subset of the product $\prod_{i\ge1}K_{F_i}^*$ a $k$-{\it slice}
if it is of the form $\{ (p_1,\dots,p_k) \}\times\prod_{i\ge k+1}K_{F_i}^*$,
where $p_i\in K_{F_i}^*$ for $1\le i\le k$. Consider any metric on
the product $\prod_{i\ge1}K_{F_i}^*$ compatible with the product topology.
Then for each $\epsilon>0$ there is $k_0$ such that
each $k_0$-slice has diameter $<\epsilon$.
Observe that for $1\le i\le k_0$, each of the maps 
$f_{F_i\#}:\#\Theta\to K_{F_i}^*$ sqeezes each subset from the
collection $K_t:t\in V_T\setminus V_{F_{k_0}}$ to a point.
As a consequence, viewing 
$\lim_{\longleftarrow}{\cal S}'=\lim_{\longleftarrow}{\cal S}_\Theta$
as a subspace in $\prod_{i\ge1}K_{F_i}^*$, we realize that for
$t\in V_T\setminus V_{F_{k_0}}$ the images $\beta(K_t)$
are contained in $k_0$-slices, and thus they all have diameters $<\epsilon$
for the metric restricted from the product.
It follows that the family $\beta(K_t):t\in V_t$ is null in
$\lim_{\longleftarrow}{\cal S}'$, and hence the family
$K_t:t\in V_T$ is null in $\lim\Theta$.

To see that the family $\Sigma_e:e\in O_T$ is null, we apply the following
general observation: {\it if a family of subsets $A_i:i\ge1$ of a metric space
is null, and if for each $i\ge1$ the family $B_{i,j}:j\ge1$ of subsets
of $A_i$ is null, then the full family $B_{i,j}:i\ge1,j\ge1$ is also null.}
We omit further details, thus completing the proof.

\medskip\noindent
{\bf 1.D.2 Example: limit of a tree system of spheres.}
For arbitrary $n\ge1$, consider a {\it tree system of $n$-spheres},
i.e. a tree system $\cal M$ of manifolds, as described in 1.B.2,
for which all manifolds $M_t:t\in V_T$ are the $n$-spheres $S^n$.
Next result is an application of Proposition 1.D.1.

\medskip\noindent
{\bf 1.D.2.1 Lemma.}
{\it The limit of any tree system of $n$-spheres is homeomorphic to $S^n$.}

\medskip
To prove Lemma 1.D.2.1, we need a special case of the 
following well known result (which we will also need later, in its full generality).

\medskip\noindent
{\bf 1.D.2.2 Lemma.}
{\it Let $M$ be an $n$-dimensional compact 
topological manifold with boundary,
and let $\cal D$ be a null and dense family of pairwise disjoint
collared $n$-disks contained in the interior of $M$.
Let $M/{\cal D}$ be the quotient space obtained by collapsing
all disks $D\in{}\cal D$ to points, i.e. the quotient space
of the decomposition of $M$ induced by $\cal D$.
Then $M/{\cal D}$ is homeomorphic to $M$, via a homeomorphism
which is identical on $\partial M$.}

\medskip\noindent
{\bf Proof:}
It follows from a theorem of Bing (Theorem 7.2 in [Fr]) that the
decomposition of $M$ induced by $\cal D$ is shirinkable
(see Section II.5 of [Dav] for the definition of shrinkability).
By Theorem 5.3 in [Dav], this implies that the quotient map
$M\to M/{\cal D}$ can be approximated by homeomorphisms,
which clearly implies our assertion. 

\medskip\noindent
{\bf Proof of Lemma 1.D.2.1:}
Let ${\cal M}=(T,\{K_t\},\{ \Sigma_e \},\{ \phi_e \})$ be any tree
system of $n$-spheres.
Note that, due to Lemma 1.D.2.2, 
each space $K_F^*:F\in{\cal F}_T$
is homeomorphic to $S^n$.
Moreover, for each finite subtrees $F\subset F'$ of $T$ the map
$f_{FF'}:K_{F'}^*\to K_F^*$ is the quotient map from $S^n$
to $S^n$ modulo a decomposition induced by a finite collection of closed collared $n$-disks. Since such a map is obviously a near homeomorphism,
we get that the associated standard inverse system ${\cal S}_{\cal M}$
consists of $n$-spheres and near homeomorphisms.

Let $(F_i)_i\ge1$ be an increasing sequence of finite subtrees of $T$
such that $\bigcup_{i\ge1}F_i=T$. Since this sequence forms a cofinal
subposet of the poset ${\cal F}_T$, the inverse sequence ${\cal S}'$
obtained from ${\cal S}_{\cal M}$ by restriction to this sequence 
has the same inverse limit. Since ${\cal S}'$ 
consists of $n$-spheres and near homeomorphisms,
it follows from a result of M. Brown (Theorem 4 in [Br]) that $\lim_{\longleftarrow}{\cal S}'\cong S^n$.
Applying Proposition 1.D.1, we get
$$
\lim{\cal M}\cong\lim_{\longleftarrow}{\cal S}_{\cal M}=
\lim_{\longleftarrow}{\cal S}'\cong S^n,
$$
as required.

\medskip\noindent
{\bf 1.D.3 Estimates of the dimension of $\lim\Theta$.}
In the remaining part of this section we use Proposition 1.D.1 to 
estimate the topological dimension of the limit
of a tree system. Some further estimates (or rather exact calculations)
will be provided in Section 2.D.

An obvious estimate, implied by the fact that each
constituent space $K_t$ embeds in $\lim\Theta$, is
$$
\dim(\lim\Theta)\ge\sup\{\dim(K_t):t\in V_T  \}.  \leqno{(1.\hbox{D}.3)}
$$
Below we provide an upper bound for $\dim(\lim\Theta)$. 
Clearly, if $\sup\{\dim(K_t):t\in V_T  \}=\infty$
then $\dim(\lim\Theta)=\infty$ as well. Thus, we restrict
our attention to tree systems which have a universal finite upper bound
for the dimensions of the constituent spaces $K_t$.

\medskip\noindent
{\bf 1.D.3.1 Proposition.}
{\it Let $\Theta=(T,\{ K_t \},\{\Sigma_e \},\{ \phi_e \})$
be a tree system of metric compacta such that 
$$
\sup\{\dim(K_t):t\in V_T  \}=n<\infty.
$$
Then $\dim(\lim\Theta)\le n+1$.}

\medskip\noindent
{\bf Proof:}
By the closed sum theorem in dimension theory, for
any $F\in{\cal F}_T$ we have $\dim(K_F)\le n$, and hence also
$\dim(K_F\setminus\bigcup_{e\in N_t}\Sigma_e)\le n$.
Recall that, for each $F\in{\cal F}_T$, we denote by $P_F$
the set of all points in $K_F^*$ which are obtained by shrinking
the subsets $\Sigma_e:e\in N_F$. Since $K_F^*\setminus P_F$
is homoemorphic to $K_F\setminus\bigcup_{e\in N_t}\Sigma_e$,
we get $\dim(K_F^*\setminus P_F)\le n$.
Since each $P_F$ is countable, it has dimension $\le0$,
and by the addition theorem we get $\dim(K_F^*)\le n+1$.
By Proposition 1.D.1, and by the properties of inverse limits,
we get $\dim(\lim\Theta)=\dim(\lim_{\longleftarrow}{\cal S}_\Theta)
\le\sup_{F\in{\cal F}_T}\dim(K_F^*)\le n+1$, as required.

\medskip\noindent
{\bf 1.D.3.2 Corollary.}
{\it Let $\Theta=(T,\{ K_t \},\{\Sigma_e \},\{ \phi_e \})$
be a tree system of metric compacta such that 
$$
\sup\{\dim(K_t):t\in V_T  \}=n<\infty.
$$
Then $\dim(\lim\Theta)\in \{ n,n+1 \}$.}

\medskip
In Section 2.D we show that both asserted in the above corollary 
values for the dimension appear, for various classes of examples.

\bigskip\noindent
{\bf 1.E An isomorphism of tree systems of spaces.}

\medskip
Let $\Theta=(T,\{ K_t \},\{\Sigma_e\},\{\phi_e\})$ and 
$\Theta'=(T',\{ K'_t \},\{\Sigma'_e\},\{\phi'_e\})$ be two tree systems
of spaces. An {\it isomorphism} $F:\Theta\to\Theta'$  is a tuple
$F=(\lambda, \{ f_t \})$ such that:
\item{(I1)} $\lambda:T\to T'$ is an isomorphism of trees;
\item{(I2)} for each $t\in V_T$ the map $f_t:K_t\to K'_{\lambda(t)}$ is a homeomorphism;
\item{(I3)} for each $e\in N_t$ we have 
$f_t(\Sigma_e)=\Sigma'_{\lambda(e)}$;
\item{(I4)} for each $e\in N_t$  the following commutation rule holds: 
$\phi'_{\lambda(e)}\circ(f_{\alpha(e)}|_{\Sigma_{\alpha(e)}})=f_{\omega(e)}\circ\phi_e$.

\medskip
An easy consequence of the definition of the limit (of a tree system of metric
compacta) is the following.

\medskip\noindent
{\bf 1.E.1 Lemma.} {\it If $\Theta,\Theta'$ are isomorphic tree systems
of metric compacta then their limits $\lim\Theta$ and $\lim\Theta'$
are homeomorphic.}

\bigskip\noindent
{\bf 1.E.2 Toru\'nczyk's Lemma and dense tree systems of oriented manifolds.}

The following result proved by Henryk Toru\'nczyk (see [J1]) has interesting
consequences concerning existence of isomorphisms for certain natural classes
of tree systems of manifolds. (Throughout this paper, by a manifold
we mean a topological manifold.)
Recall that a family of subsets of a topological space is {\it dense}
if the union of this family is a dense subset.

\medskip\noindent
{\bf 1.E.2.1 Toru\'nczyk's Lemma.}
{\it Let $M$ be a compact $n$-dimensional topological manifold with or without boundary,
and let ${\cal D},{\cal D}'$ be two families of collared $n$-disks
in $\hbox{int}(M)$ such that each family consists of pairwise disjoint sets
and both families are null and dense. Then each homeomorphism
$h:\partial M\to\partial M$ which is extendable to a homeomorphism of $M$
admits an extension to a homeomorphism $H:M\to M$ which maps ${\cal D}$
to ${\cal D}'$. More precisely, this means that there is an associated
bijective map $\nu:{\cal D}\to{\cal D}'$ such that for each
$\Delta\in{\cal D}$ the restriction $H|_\Delta$ maps $\Delta$
homeomorphically on the disk $\nu(\Delta)\in{\cal D}'$.}

\medskip
The above lemma provides motivation for the following.

\medskip\noindent
{\bf 1.E.2.2 Definition.}
Let $\cal M$ be a tree system of closed manifolds, with families
of manifolds $\{M_t\}$ and disks $\{\Delta_e\}$ as in Example 1.B.2.
We say that this system is {\it dense} if for each $t\in V_T$
the family ${\cal D}_t=\{ \Delta_e:e\in N_t \}$ is dense in the
manifold $M_t$.

\medskip
We denote the constituent spaces of a dense
system of manifolds $M_t$ by $M_t^\circ$. The symbol 
$M_t^\circ$ is meant to contain information both of the space itself
and of the peripheral subspaces contained in
this space. Given $M_t$, it follows from Lemma 1.E.2.1 
that the corresponding $M_t^\circ$ is unique up to
a homeomorphism preserving the peripheral subspaces.
We will call any constituent space of form $M_t^\circ$
(viewed again as equipped with its standard family of peripheral
subspaces) a {\it densely punctured manifold $M_t$}.

\medskip
A recursive application of Toru\'nczyk's Lemma, together with
Lemma 1.E.1, immediately yield the following.

\medskip\noindent
{\bf 1.E.2.3 Proposition.} 
{\it Let $M$ be a closed connected oriented topological manifold. Let $\cal M$, 
${\cal M}'$ be two dense tree systems of manifolds such that} 
\item{(1)} {\it all manifolds
in the corresponding families $\{ M_t \}$ and $\{ M'_{t'} \}$ are
homeomorphic to $M$,} 
\item{(2)} {\it all  maps 
$\phi_e:\Sigma_e=\partial\Delta_e\to
\Sigma_{\bar e}=\partial\Delta_{\bar e}$ respect orientations,
i.e. reverse the induced orientations on the corresponding spheres,
and the same holds for all maps $\phi'_{e'}$.}

\noindent
{\it Then the tree systems $\cal M$, ${\cal M}'$ are isomorphic,
and consequently, their limits are homeomorphic.}

\medskip
Note that,
due to the above proposition, 
each closed connected oriented manifold $M$
determines uniquely up to isomorphism the dense tree system of manifolds
satisfying conditions (1) and (2) as in the proposition. 
We denote this tree system by ${\cal M}(M)$ and call it
{\it the dense tree system of manifolds $M$} or
{\it the Jakobsche system for $M$}. 
The latter term is motiveted by the fact that the tree systems
${\cal M}(M)$ are intimately related to some inverse sequences
described by W. Jakobsche in [J2]. We describe this relationship in Part 2 of the paper, especially in Example 2.C.7.

As a consequence
of Lemma 1.E.1, $M$ as above determines, uniquely up to homeomorphism,
the compact metric space $\lim{\cal M}(M)$, which we denote by
${\cal X}(M)$ and call {\it the tree of manifolds $M$} or 
{\it the Jakobsche space for $M$}. In Jakobsche's paper [J2]
this space is denoted by $X(M,\{M\})$, and it is obtained as
inverse limit of an inverse sequence mentioned in the previous paragraph.

\medskip\noindent
{\bf Remark.} 
In [J2] Jakobsche considered also a more general class of 
spaces obtained as inverse limits and
determined uniquely up to homeomorphism by a finite or infinite family
$\cal N$ of closed connected oriented manifolds of the same dimension.
We discuss the corresponding dense tree systems of manifolds, 
for finite families $\cal N$, 
in Example 3.A.2.

\medskip\noindent
{\bf 1.E.3 Example: tree of spheres is a sphere.}

For arbitrary $n\ge1$, consider the dense tree system 
${\cal M}(S^n)$ of $n$-spheres, as described in Subsection 1.E.2.
Note that this tree system can be alternatively thought of as follows.

\medskip\noindent
{\bf 1.E.3.1 Remark: an alternative description of the tree system 
${\cal M}(S^n)$.}

Recall that the unique topological space $(S^n)^\circ$ obtained from $S^n$
by deleting interiors of  $n$-disks $D$ from any null and dense
family $\cal D$ consisting of pairwise disjoint collared disks
is called the $(n-1)$-dimensional Sierpi\'nski compactum. (The
uniqueness follows from [Ca] if $n\ne4$, and in the remaining case the argument in [Ca] also holds true in view of the later proofs of the Annulus Theorem and the Approximation Theorem for $n=4$, due to F. Quinn [Qu];
in fact, this follows also from Toru\'nczyk's Lemma 1.E.2.1.) 
The space $(S^n)^\circ$ contains a family of distinguished subsets, called
{\it peripheral spheres}, which coincides with
the family of boundaries $\partial D$ of the disks $D\in{\cal D}$.
The tree system ${\cal M}(S^n)$ can be described as follows.
It is the unique tree system of metric compacta in which
all constituent spaces are $(n-1)$-dimensional Sierpi\'nski compacta,
and families of peripheral subsets in all of these spaces coincide with
the families of peripheral spheres.

\medskip
Next result is a special case of Lemma 1.D.2.1.

\medskip\noindent
{\bf 1.E.3.2 Corollary.}
{\it The limit $\lim {\cal M}(S^n)$, i.e. the tree of spheres
${\cal X}(S^n)$, is homeomorphic to $S^n$.}

\medskip
A different argument for proving Corollary 1.E.3.2 is sketched in
Remark 3.B.12.2.

\bigskip\noindent
{\bf 1.E.4 Example: the tree of non-orientable manifolds $N$.}

Reall the following rather well known fact.

\medskip\noindent
{\bf 1.E.4.1 Lemma.}
{\it Let $N$ be a closed connected non-orientable 
topological manifold of
dimension $n$
and let $D,D'$ be any collared $n$-disks in $N$.
Then each homeomorphism $D\to D'$ (no matter how it behaves with
respect to local orientations in the disks) extends to a homeomorphism
of $N$. }

\medskip
A recursive application of Toru\'nczyk's Lemma 1.E.2.1 and
Lemma 1.E.4.1 yields the following.

\medskip\noindent
{\bf 1.E.4.2 Proposition.}
{\it Let $N$ be a closed connected non-orientable topological manifold,
and let $\cal M$, ${\cal M}'$ be any two dense tree systems of manifolds 
such that all manifolds
in the corresponding families $\{ M_t \}$ and $\{ M'_{t'} \}$ are
homeomorphic to $N$. 
Then the tree systems $\cal M$, ${\cal M}'$ are isomorphic,
and consequently, their limits are homeomorphic.}

\medskip
The tree system as above, uniquely determined by $N$,
will be denoted ${\cal M}(N)$ and called the 
{\it dense tree system of manifolds $N$}, or the
{\it Jakobsche system for $N$}. Its limit, denoted ${\cal X}(N)$,
will be called the {\it tree of manifolds $N$} or the
{\it Jakobsche space for $N$}.

\medskip\noindent
{\bf 1.E.4.3 Remark.}
Proposition 1.E.4.2 clarifies the picture with trees of non-orientable
manifolds, as documented so far in the literature. Namely, in the paper [St]
Paul Stallings describes these spaces in a way which could be translated
to our setting as follows: the tree of manifolds $N$ is the limit of
a dense tree system as in Proposition 1.E.4.2 satisfying some
additional technical condition for the connecting maps called
{\it dense orientation condition} (we do not recall this condition).
Proposition 1.E.4.2 shows that this additional condition plays
in fact no role. Moreover, Proposition 1.E.4.2 applies to all topological manifolds, while the methods used by Stallings in [St] allowed him to deal only with PL manifolds.


\bigskip\bigskip
\centerline{\bf 2. Extensions of tree systems and 
associated inverse systems.}

\bigskip
In this part of the paper we present some alternative expressions of
trees of spaces as inverse limits. It turns out that these expressions
(rather than those related to standard associated inverse systems
described in Section 1.D)
allow to relate the spaces obtained as limits of
tree systems with certain
previously studied classes of topological spaces defined in terms
of inverse limits (e.g. Jakobsche trees of manifolds, Markov
compacta as defined in [Dr], etc.). 
Moreover, these expressions seem to be potentially more convenient
for the purpose of recognizing ideal boundaries of various
spaces and groups as (homeomorphic to) some specific trees of spaces
(see Remark 2.C.8 for 
examples of such applications occuring in the literature).

The alternative expressions of trees of spaces as inverse limits, 
presented in this section,
turn out to be useful for exact calculations of topological dimension
for certain classes of trees of spaces (see Section 2.D).

\bigskip\noindent
{\bf 2.A Extended spaces and maps.}

\medskip
Let $\Theta=(T,\{ K_t \},\{ \Sigma_e \},\{ \phi_e \})$
be a tree system of metric compacta.
Suppose that for each $e\in O_T$
we are given a compact metric space $\Delta_e$, its compact subspace
$S_e$, and a homeomorphism $\varphi_e:S_e\to\Sigma_e$.

For any $t\in V_T$ let $T(t)$ be a subtree of $T$ spanned
on the set $\{t\}\cup\{ \omega(e):e\in N_t \}$ (i.e. $t$ and all vertices
adjacent to $t$).  We define a tree system 
$$
\Theta(t)=(T(t), \{ K'_s \}, \{ \Sigma'_e \}, \{ \phi'_e \})
$$
as follows. Put $K'_t=K_t$ and $K'_{\omega(e)}=\Delta_e$ for each 
$e\in N_t$. Put $\Sigma'_e=\Sigma_e$ and $\Sigma'_{\bar e}=S_e$
for each $e\in N_t$. Finally, put $\phi'_e=\varphi_e^{-1}$
and $\phi'_{\bar e}=\varphi_e$ for each $e\in N_t$.
Denote by $\hat K_t=\lim\Theta(t)$ the limit of the above tree system.
Observe that $K_t$ and the sets $\Delta_e:e\in N_t$ are then canonically 
the subspaces of $\hat K_t$ (in particular, they are compact subspaces).
Moreover, the sets $\Sigma_e$ and $S_e$, viewed in the above way as subspaces
of $\hat K_{t}$, do coincide. 
We call any family $\hat K_t:t\in V_t$
as above a {\it family of extended spaces} for $\Theta$.

\medskip\noindent
{\bf 2.A.1 Examples.}

\item{(1)} For each $e\in O_T$ put $\Delta_e=S_e=\Sigma_e$
and $\varphi_e=id_{\Sigma_e}$. This defines what we call the {\it trivial
family of extended spaces} for $\Theta$. Note that in this family we have
$\hat K_t=K_t$ for each $t\in V_T$.

\item{(2)} Let $\cal M$ be a tree system of manifolds as in Example 1.B.2.
Let $\Delta_e:e\in O_T$ be the family of $n$-disks as in this example,
and put $S_e=\partial\Delta_e$  and $\varphi_e=id_{S_e}$ for each such $e$.
Observe that then we have
$\hat K_t=M_t$  for each $t\in V_T$.
Indeed, by Lemma 1.D.2.2 all reduced partial unions of the system
$\Theta(t)$ corresponding to finite subtrees $F\subset T(t)$ containing $t$
(i.e. all spaces in the standard associated inverse system 
${\cal S}_{\Theta(t)}$ restricted to the cofinal subposet in
${\cal F}_{T(t)}$ consisting of all subtrees containing $t$) 
are then homeomorphic to $M_t$. 
Moreover, an argument as
in the proof of Lemma 1.D.2.2 shows that the maps in 
the above restriction of the system
${\cal S}_{\Theta(t)}$
are all near homeomorphisms, and thus (applying also Proposition 1.D.1)
we get $\hat K_t=\lim\Theta(t)\cong\lim_{\longleftarrow}{\cal S}_{\Theta(t)}
\cong M_t$.
(A more canonical identification of $\hat K_t$ with $M_t$
is easily provided by Theorem 3.B.10.)
This means that as a family of extended spaces for $\cal M$
we can take the initial family of manifolds $M_t$.
We will call the family $\{ \hat K_t=M_t \}$ as above the
{\it standard family of extended spaces for $\cal M$.}

\item{(3)} The previous example can be generalized as follows.
Let $\Theta=(T,\{ K_t \},\{ \Sigma_e \},\{ \phi_e \})$
be any tree system of metric compacta. For each $e\in N_T$
put $\Delta_e=cone(\Sigma_e)$ and $S_e$ to be the base
of this cone, and let $\varphi_e$ be the identity on $S_e=\Sigma_e$.
We call the associated family $\hat K_t:t\in V_T$ {\it the family
of conically extended spaces} for $\Theta$.

\medskip
Given a family of extended spaces for $\Theta$ determined,
as above, by a family of pairs $(\Delta_e,S_e):e\in O_T$, 
consider the family $\hat K_e:e\in O_T$ of subspaces defined by 
$\hat K_e=\hat K_{\omega(e)}\setminus(\Delta_{\bar e}\setminus S_{\bar e})$.
Suppose
that we are given a family of maps $\delta_e:\hat K_e\to\Delta_e$
such that $\delta_e|_{\Sigma_{\bar e}}=\phi_{\bar e}$.
Any tuple ${\cal E}=(\{ \hat K_t:t\in V_T \},\{ \delta_e:e\in O_T \})$
as above will be called {\it a family of extended spaces
and maps} for $\Theta$.

To any family $\cal E$ of spaces and maps as above
we associate an inverse system of metric compacta, denoted ${\cal S}_{\cal E}$,
as follows.
Let $({\cal F}_T,\subset)$ be the poset of all finite subtrees of $T$.
For any $F\in{\cal F}_T$ let $\hat K_F$ be the quotient of the disjoint union
$$
\hat K_F=\bigcup_{t\in V_F}\big[\hat K_t\setminus\bigcup_{e\in N_t\cap O_F}
(\Delta_e\setminus S_e)\big]/\sim
$$
where $\sim$ is the equivalence relation induced by the equivalences
of form $x\sim\phi_e(x)$ for all $e\in O_F$
and all $x\in\Sigma_e$ (where we view $\Sigma_e$
and $\Sigma_{\bar e}$ canonically as subsets in 
$\hat K_{\alpha(e)}$ and $\hat K_{\omega(e)}$, respectively).
This gives us the family $\hat K_F:F\in{\cal F}_T$ of spaces
in the inverse system ${\cal S}_{\cal E}$.

We now turn to defining the maps $h_{F_1F_2}:\hat K_{F_2}\to \hat K_{F_1}$,
for all pairs $F_1\subset F_2$ of finite subtrees of $T$. 
We do this in two steps,
as follows. First, suppose that $V_{F_2}=V_{F_1}\cup\{ t \}$
for some $t\notin V_{F_1}$.
Let $e$ be the unique oriented edge in $F_2$ with $\omega(e)=t$.
Viewing canonically $\Delta_e$ as a subset of $\hat K_{F_1}$, 
we put
$$
h_{F_1F_2}(x)=\cases{x & for $x\in \hat K_{F_1}
                                 \setminus(\Delta_e\setminus S_e)
                                  \subset\hat K_{F_2}$\cr
                                  \delta_e(x) & for $x\in\hat K_e\subset\hat K_{F_2}$}.
$$
In the second step, 
for any pair $F'\subset F$ of finite subtrees of $T$ we consider
a sequence $F_1\subset F_2\subset\dots\subset F_m$ such that
$F_1=F', F_m=F$ and the pairs $F_i,F_{i+1}$ are of the form
as in the first step. We then put
$$
h_{FF'}=h_{F_1F_2}\circ h_{F_2F_3}\circ\dots\circ h_{F_{m-1}F_m}
$$
and we observe that the resulting map does not depend on the choice
of the above sequence $F_1,\dots,F_m$ (which is in general not unique).
A related observation is that $h_{F_1F_3}=h_{F_1F_2}\circ h_{F_2F_3}$
whenever $F_1\subset F_2\subset F_3$.

\medskip\noindent
{\bf 2.A.2 Definition.}
Let $\Theta$ be a tree system of metric compacta, and let $\cal E$
be an associated family of extended spaces and maps for $\Theta$.
The {\it associated inverse system for $\Theta$  induced by $\cal E$}
is the tuple 
$$
{\cal S}_{\cal E}=(\{ \hat K_F:F\in{\cal F}_T \},\{ h_{FF'}:F\subset F' \}).
$$

\medskip
In the next section we show (see Theorem 2.B.4) that in some cases
one can use limit of an associated inverse system ${\cal S}_{\cal E}$
as an alternative description of the limit of a tree system.

\medskip\noindent
{\bf 2.A.3 Remark.}
Any associated inverse system ${\cal S}_{\cal E}$ as above has
various cofinal subsequences. 
Since any cofinal subsystem
is viewed as equivalent to the original one (for example, its inverse limit
is canonically the same), this allows simpler (in certain sense) descriptions
of the system ${\cal S}_{\cal E}$, 
which might be more convenient for some purposes.
More precisely, let $F_1\subset F_2\subset\dots$ be an increasing
sequence of finite subtrees of $T$. Clearly, this sequence is cofinal with
${\cal F}_T$ iff $\cup_{i=1}^\infty V_{F_i}=V_T$. For any such cofinal sequence
$(F_i)$ the tuple 
$$
{\cal S}_{{\cal E},(F_i)}=(\{ \hat K_{F_i}:i\in N\hskip-3pt I \},\{ h_{F_jF_i}:i<j \})
$$
is an inverse sequence equivalent to the inverse system ${\cal S}_{\cal E}$.

\bigskip\noindent
{\bf 2.B. Relationship to the limit of a tree system.}

\medskip
In order to relate the limit $\lim\Theta$ with the inverse limit of some
associated inverse system ${\cal S}_{\cal E}$, we need one more
condition on the associated family $\cal E$ of extended spaces and maps.
We use the notation as in the previous section.

Let $\gamma=(t_0,t_1,\dots,t_m)$ be any finite combinatorial path in $T$
of length $m\ge2$.
For $i=1,\dots,m$ let $e_i=(t_{i-1},t_i)$ be the consecutive oriented edges
in $\gamma$. We then have
$$
\Delta_{e_m}\subset{\hat K_{e_{m-1}}}\,\,\,\,_{\longrightarrow}
\hskip-18pt^{\delta_{e_{m-1}}}\,\,\Delta_{e_{m-1}}\subset\,\,\dots
\,\,\,_{\longrightarrow}\hskip-14pt^{\delta_{e_2}}\,\,\,\,
\Delta_{e_2}\subset\hat K_{e_1}\,_{\longrightarrow}\hskip-14pt
^{\delta_{e_1}}\,\,\,\,\Delta_{e_1}
$$
and we denote by $\delta_\gamma:\Delta_{e_m}\to\Delta_{e_1}$
the composition map $\delta_\gamma=\delta_{e_1}\circ\dots\circ
\delta_{e_{m-2}}\circ\delta_{e_{m-1}}|_{\Delta_{e_m}}$.

\medskip\noindent
{\bf 2.B.1 Definition.}
We say that an associated family $\cal E$ of extended spaces and maps is 
{\it fine} if for each $e\in O_T$ the family of images 
$\delta_\gamma(\Delta_{e_m})$ (where $e_m$ is the terminal edge
in $\gamma$), for all combinatorial paths $\gamma$ in $T$
of length $\ge2$ and starting with $e_1=e$, is a null family of subsets
in $\Delta_e$.

\medskip\noindent
{\bf Remark.} The condition of fineness does not follow automatically from
the nullity conditions in $\Theta$ (for families $\{ \Sigma_e:e\in N_t \}$).

\medskip
Next definition describes a subclass of fine families $\cal E$
of extended spaces and maps which occur in practical situations
where we express limits of tree systems as inverse limits.

\medskip\noindent
{\bf 2.B.2 Definition.}
Let ${\cal E}=(\{ \hat K_t \},\{ \delta_e \})$ be a family of extended 
spaces and maps related to a family of pairs $\{ (\Delta_e,S_e) \}$.
We say that $\cal E$ is {\it contracting} if there are metrics
$\hat d_t$ on the extended spaces $\hat K_t$ and 
a constant $0\le c<1$ such that

\item{(1)}
for each $e,e'\in O_T$ such that $e'\ne\bar e$ and $\alpha(e')=\omega(e)$
the restricted map $\delta_e|_{\Delta_{e'}}$ is a $c$-contraction
with respect to the metrics $\hat d_{\omega(e)}$ in 
$\Delta_{e'}\subset\hat K_{\omega(e)}$ and $\hat d_{\alpha(e)}$
in $\Delta_e\subset\hat K_{\alpha(e)}$;
more precisely, for any $x,y\in\Delta_{e'}$ we have
$$
\hat d_{\alpha(e)}(\delta_e(x),\delta_e(y))\le c\cdot \hat d_{\omega(e)}(x,y);
$$

\item{(2)} if $c>0$ then there is another constant $C>0$ such that
for each $t\in V_T$ we have
$\hbox{diam}(\hat K_t,\hat d_t)\le C$.

\medskip
We omit the proof of the following apparent fact.

\medskip\noindent
{\bf 2.B.3 Fact.} 
{\it Any contracting family $\cal E$ of extended spaces and maps
is fine.}

\medskip\noindent
{\bf 2.B.4 Theorem.}
{\it Let $\Theta$ be a tree system of metric compacta, let $\cal E$
be an associated family of extended spaces and maps for $\Theta$,
and suppose that $\cal E$ is fine. Let ${\cal S}_{\cal E}$ 
be the associated inverse system
of compact topological spaces. Then the limit $\lim\Theta$ is canonically
homeomorphic to the inverse limit $\lim_\leftarrow {\cal S}_{\cal E}$.}

\bigskip\noindent
{\bf Proof:} 
The proof consists of three parts. In the first part we describe explicitely
strings $(x_F)_{F\in{\cal F}_T}$ which represent points of
$\lim_\leftarrow {\cal S}_{\cal E}$. In the second part, we use this decription
to define a natural map $\beta:\lim_\leftarrow {\cal S}_{\cal E}\to\lim\Theta$,
which is a bijection. Finally, in the third part we prove that
$\beta$ is a homeomorphism.

Recall that, by definition of inverse limit, each element
$x\in\lim_\leftarrow {\cal S}_{\cal E}$ is a tuple $(x_F)_{F\in{\cal F}_T}$,
with $x_F\in\hat K_F$, such that $h_{FF'}(x_F)=x_{F'}$ whenever $F'\subset F$.
We will show that there are two kinds of such tuples in
$\lim_\leftarrow {\cal S}_{\cal E}$.
At first, we consider tuples $(x_F)$ satisfying the following property:

\smallskip
\item{($\star$)} for any cofinal subsequence 
$F_1\subset F_2\subset\dots$ in the poset ${\cal F}_T$
the sequence $x_{F_n}$ {\it eventually stabilizes},
i.e. there is $N$ such that for all $n\ge N$ we have 
$x_{F_n}\in K_{F_n}\subset\hat K_{F_n}$
and, under the natural inclusions $K_{F_n}\subset K_{F_{n+1}}$, 
we have the equalities $x_{F_{n+1}}=x_{F_n}$.

\smallskip
\noindent
A large class of such tuples can be described as follows.
Let $y\in K_t$ for some $t\in V_T$. This $y$ determines the tuple
$(x^y_F)$ as follows. If $F$ contains $t$, we have the canonical inclusion
$K_t\subset\hat K_F$, and we put $x^y_F=x$.
If $F$ does not contain $t$, let $F'$ be the smallest subtree of $T$
containing both $F$ and $t$. Viewing $K_t$ as a subset of
$\hat K_{F'}$ we put $x^y_F:=h_{F'F}(x)$. It is easy to check
that $x^y=(x^y_F)$ is then a string, 
i.e. $x^y\in\lim_\leftarrow {\cal S}_{\cal E}$, and it obviously satisfies
condition $(\star)$. Moreover, we have the following.

\medskip\noindent
{\bf Claim 1.} 
{\it Each string $(x_F)\in\lim_\leftarrow {\cal S}_{\cal E}$ which satisfies property ($\star$) has form $(x^y_F)$, as described above.}

\medskip
The claim follows by observing that if $y$ is the element to
which some sequence $x_{F_n}$ stabilizes, then any other
such sequence stabilizes to the same element, and consequently
the string necessarily has the asserted form $(x^y_F)$.

We now turn to strings $(x_F)$ not satisfying condition $(\star)$.
A large class of such strings is induced by ends $z\in E_T$.
Given any $z\in E_T$,
we define the tuple $(x^z_F)$ as follows.
For any $F\in{\cal F}_T$ let $\rho_F$ be the minimal (with respect to inclusion)
ray starting at a vertex of $F$ and such that $[\rho_F]=z$.
Let $\gamma_n=(e_1',\dots,e_n')$ be the initial path of length $n$ 
in $\rho_F$, and let 
$\delta_{\gamma_n}:\Delta_{e_n'}\to\Delta_{e_1'}$ 
be the map described at the beginning of this section.
Note that the sequence $(Q_n)$ of compact 
subsets in $\Delta_{e_1'}$
given by $Q_n=\delta_{\gamma_n}(\Delta_{e_n'})$ is then
decreasing and, due to fineness of $\cal E$, we have
$\lim_{n\to\infty}\hbox{diam}(Q_n)=0$. Consequently,
the intersection $\cap Q_n$ is a single point
in $\Delta_{e_1'}\subset\hat K_F$, and we take it as $x^z_F$.
We make the following observations (omitting their straightforward proofs):

\item{(1)} $(x^z_F)$ is a string, i.e. it belongs to
$\lim_\leftarrow {\cal S}_{\cal E}$,
\item{(2)} $(x^z_F)$ does not satisfy condition $(\star)$,
\item{(3)} if $z\ne z'$ then $(x^z_F)\ne(x^{z'}_F)$.

\medskip\noindent
{\bf Claim 2.} 
{\it Each element $(x_F)\in\lim_\leftarrow {\cal S}_{\cal E}$
which does not satisfy condition $(\star)$ has a form
$(x^z_F)$ as above, for some end $z\in E_T$.
Moreover, $z\to(x^z_F)$ is a bijective correspondence
between the set of ends of $T$ and the set of elements of 
$\lim_\leftarrow {\cal S}_{\cal E}$ not satisfying $(\star)$.}

\medskip
To prove the claim, note that
for each string $(x_F)$ not satisfying $(\star)$ there is a cofinal sequence 
$F_1\subset F_2\subset\dots$ for which elements $x_{F_n}$
change infinitely often. By passing to a subsequence,
we may assume that $x_{F_{n+1}}\ne x_{F_n}$ 
(i.e. it is not true that 
$x_{F_{n+1}}=x_{F_n}\in K_{F_n}\subset K_{F_{n+1}}$) for all $n$.
Observe that in this situation we have $x_{F_n}\in\Delta_{e_n}$
for some unique $e_n\in N_{F_n}$, and that the edges $e_n$
induce uniquely a ray $\rho$ in $T$ of form
$\rho=(e_1,\dots,e_2,\dots,e_3,\dots)$. 
It is not hard to realize that then $(x_F)=(x^{[\rho]}_F)$,
which yields the first assertion of the claim.
The second assertion follows then from observation (3) stated
just before Claim 2.

\medskip
We are now ready to describe a natural map 
$\beta:\lim_\leftarrow {\cal S}_{\cal E}\to\lim\Theta$, which will be
our candidate for a homeomorphism.
If a string $x=(x_F)$ satisfies condition $(\star)$,
there is a point $y$ to which $(x_F)$ stabilizes,
and this $y$ belongs to some $K_t$. Viewing $K_t$ as a subset in
$\lim\Theta$, we put $\beta(x)=y$. If $x$ does not satisfy $(\star)$,
it corresponds uniquely to some end $z$ of $T$. Viewing ends of $T$
as elements of $\lim\Theta$, we put $\beta(x)=z$.
We skip a direct verification of the fact that so described $\beta$
is well defined, and that it is a bijection.

Since both spaces $\lim_\leftarrow {\cal S}_{\cal E}$ and $\lim\Theta$
are compact, to prove that $\beta$ is a homeomorphism,
it is sufficient to show that it is continuous. To do this, we will show that
for any open set $G(U)\in{\cal B}$ its preimage 
$\beta^{-1}(G(U))$ is open in $\lim_\leftarrow {\cal S}_{\cal E}$.
More precisely, viewing $\lim_\leftarrow {\cal S}_{\cal E}$ as subspace
of the product $\prod_F\hat K_F$, we will show that
$\beta^{-1}(G(U))=W\cap\lim_\leftarrow {\cal S}_{\cal E}$
for some open set $W$ in $\prod_F\hat K_F$.
Suppose that $U$ is an open subset of $K_{F_0}$ saturated
with respect to the family ${\cal A}_{F_0}$. Put
$$
\hat U=U\cup\bigcup\{ \Delta_e:e\in N_{F_0}\hbox{ and }
\Sigma_e\subset U \}
$$
and note that $\hat U$ is an open subset of $\hat K_{F_0}$.
Take $W=\prod_F W_F$, where $W_{F_0}=\hat U$
and $W_F=\hat K_F$ for $F\ne F_0$. Clearly, $W$ is open in
$\prod_F\hat K_F$. It is also not hard to deduce from Claims 1 and 2,
and from the descriptions of strings,
that $\beta^{-1}(G(U))=W\cap\lim_\leftarrow {\cal S}_{\cal E}$.
This completes the proof.

\bigskip\noindent
{\bf 2.C Tree systems with ANR peripheral subspaces.}

\medskip

This section is devoted to showing that any tree system $\Theta$ in which
all peripheral subspaces $\Sigma_e$ are ANRs admits a fine
family ${\cal E}=(\{ \hat K_t \},\{ \delta_e \})$ of extended
spaces and maps in which $\{ \hat K_t \}$ is the family of conically
extended spaces (see Example 2.A.1(3)) and $\{ \delta_e \}$
is an associated family of $0$-contracting maps. See Proposition 2.C.2
for precise statement of this main result of the section.

Recall that a compact metric space $\Sigma$ is an ANR 
(absolute neighbourhood retract) if for any embedding of $\Sigma$
in another compact metric space $K$ there is a neighbourhood
of $\Sigma$ in $K$ which retracts on $\Sigma$. 
Recall also that every compact polyhedron is an ANR.
Lemma 2.C.3 below
presents a different characterization of ANR spaces, more convenient
for our purposes.

\medskip\noindent
{\bf 2.C.1 Definition.}
Let ${\cal E}=(\{ \hat K_t \},\{ \delta_e \})$ be a family of extended
spaces and maps for a tree system $\Theta$.
\item{(1)} We say that $\Theta$ is {\it peripherally ANR} if
all peripheral subspaces $\Sigma_e:e\in O_T$ of $\Theta$
are ANRs.
\item{(2)} We say that $\cal E$ is {\it conical} if the corresponding
family $\{ \hat K_t \}$ is the family of conically
extended spaces for $\Theta$.
\item{(3)} We say that $\cal E$ is {\it $0$-contracting} if
for each $e,e'\in O_T$ such that $e'\ne\bar e$ and $\alpha(e')=\omega(e)$
the restricted map $\delta_e|_{\Delta_{e'}}$ is a $0$-contraction;
equivalently, for all $e,e'$ as above $\delta_e$ maps  
the subspace $\Delta_{e'}\subset\hat K_e$ to a point.

\medskip
Note that a 0-contracting family $\cal E$ is automatically fine
(compare Fact 2.B.3) and thus, due to Theorem 2.B.4, can be used
to express the limit $\lim\Theta$ as an inverse limit.
The main result of this section is the following.

\medskip\noindent
{\bf 2.C.2 Proposition.}
{\it Each peripherally ANR tree system 
has a conical and $0$-contracting family 
of extended spaces and maps.}

\medskip
The proof of Proposition 2.C.2 requires various preparations.
We start with a lemma which characterizes ANR spaces $\Sigma$
in terms of maps to the cones over $\Sigma$. Since this characterization
is an obvious reformulation of the definition of an ANR space,
we omit its proof.

\medskip\noindent
{\bf 2.C.3 Lemma.} {\it Let $\Sigma$ be a compact metric space, and let
$cone(\Sigma)$ be the cone over $\Sigma$, with $\Sigma$ canonically
identified as the cone base. Then the following two conditions are equivalent:}
\item{(1)} {\it for any metric compactum $K$ containing $\Sigma$ 
as a subspace
there is a continuous map $f:K\to cone(\Sigma)$ such that $f|_\Sigma=id_\Sigma$;}
\item{(2)} {$\Sigma$ is an ANR.}

\medskip
We now turn to discussing decompositions of the constituent spaces
$K_t$ induced by families of their peripheral subsets.
For any $e\in O_T$, put 
$$
{\cal A}_e=\{ \Sigma_{e'}:e'\ne{\bar e}\hbox{ and } \alpha(e')=\omega(e) \}.
$$
Clearly, ${\cal A}_e$ is a null family of pairwise disjoint compact subsets
of the space $K_{\omega(e)}$. We will view ${\cal A}_e$
as a decomposition of $K_{\omega(e)}$ into closed subsets
by considering all singletons $\{x\}$ with $x\notin\cup{\cal A}_e$,
together with the sets from ${\cal A}_e$,
as elements  of this decomposition. Equivalently, we identify 
the family ${\cal A}_e$ with the decomposition of $K_{\omega(e)}$
in which ${\cal A}_e$ is the set of nondegenerate elements.
By the fact that the family ${\cal A}_e$ is null we immediately get
the following (compare [Dav], Proposition 3 on p. 14).

\medskip\noindent
{\bf 2.C.4 Fact.}
{\it For each $e\in O_T$ the decomposition ${\cal A}_e$ of the space
$K_{\omega(e)}$ is upper semicontinuous.}

\medskip
Denote by $K_{\omega(e)}/{\cal A}_e$ the quotient space
of the decomposition ${\cal A}_e$. Since any quotient of 
an upper semicontinuous decomposition of a metric space
is metrizable (see [Dav], Proposition 2 on p. 13), we get

\medskip\noindent
{\bf 2.C.5 Fact.} 
{\it For each $e\in O_T$ the quotient $K_{\omega(e)}/{\cal A}_e$
is a compact metrizable space, and the subspace $\Sigma_{\bar e}$
canonically topologically embeds in this quotient.}

\medskip
Fact 2.C.5 together with Lemma 2.C.3 yield the following.

\medskip\noindent
{\bf 2.C.6 Corollary.} 
{\it Suppose that the peripheral subset $\Sigma_{\bar e}$ is an ANR.
Then there is a continuous map 
$f_e:K_{\omega(e)}/{\cal A}_e\to cone(\Sigma_{\bar e})$ which
is identical on $\Sigma_{\bar e}$. Consequently, there is an induced
continuous map $f'_e:K_{\omega(e)}\to cone(\Sigma_{\bar e})$
which is identical on $\Sigma_{\bar e}$ and which
maps each subset $\Sigma_{e'}\in{\cal A}_e$ to a point.}

\bigskip\noindent
{\bf Proof of Proposition 2.C.2:}
Let $\{ \hat K_t:t\in V_T \}$ be the associated 
family of conically extended spaces for $\Theta$,
as described in Example 2.A.1(3),
and let $\{ \hat K_e:e\in O_T \}$ be the corresponding family
of subspces. For each $e\in O_T$ define
a map $\delta'_e:\hat K_e\to cone(\Sigma_e)$
by 
$$
\delta'_e(x)=
\cases{f'_e(x) & if $x\in K_{\omega(e)}\subset\hat K_e$\cr
f'_e(\Sigma_{e'}) & if $x\in\Delta_{e'}\subset\hat K_e$,
for some $e'\in N_{\omega(e)}\setminus\{\bar e\},$\cr}
$$
where $f'_e$ is a map as in Corollary 2.C.6.

After identifying $cone(\Sigma_{\bar e})$ with $cone(\Sigma_e)$
via cone of the map $\phi_{\bar e}:\Sigma_{\bar e}\to\Sigma_e$,
and then identifying $cone(\Sigma_e)$ with $\Delta_e$,
the above described map $\delta'_e$ gives the map
$\delta_e:\hat K_e\to\Delta_e$ such that:
\item{(1)} $\delta_e$ maps each subset 
$\Delta_{e'}\subset\hat K_e$, 
for $e'\in N_{\omega(e)}\setminus\{\bar e\}$, 
to a point;
\item{(2)} the restriction of $\delta_e$ to $\Sigma_{\bar e}$
coincides with $\phi_{\bar e}$.

\noindent
Thus, the tuple ${\cal E}=(\{ \hat K_e \},\{ \delta_e \})$
is a conical and $0$-contracting family
of extended spaces and maps for $\Theta$, as required.
This finishes the proof.

\medskip\noindent
{\bf 2.C.7 Example: Jakobsche's inverse sequences.}

As we have already mentioned in Section 1.E,
W\l odzimierz Jakobsche has introduced and studied (in his papers [J1,J2]) 
the spaces which appear in the present paper as limits of 
dense tree systems of manifolds,
in particular the tree systems of manifolds $M$, as described in 
Proposition 1.E.2.3 (the Jakobsche spaces ${\cal X}(M)$).
Jakobsche has introduced those spaces as inverse limits of certain
inverse sequences of manifolds, which we can describe in our terms
as follows.

Given a closed oriented manifold $M$, let ${\cal M}={\cal M}(M)$
be the dense tree system of manifolds $M$ (as defined in
Section 1.E, right after Proposition 1.E.2.3). 
Let ${\cal E}=(\{ \hat K_t \},\{ \delta_e \})$
be any fine conical family of extended spaces and maps for $\cal M$
(existence of which is justified by Proposition 2.C.2).
Let ${\cal S}_{\cal E}=(\{ \hat K_F \},\{ h_{F,F'} \})$
be the associated inverse system for $\cal M$,
as in Definition 2.A.2. Note that the spaces $\hat K_F$ in this system
are homeomorphic to iterated connected (oriented) sum of copies of $M$.

Jakobsche has considered cofinal inverse subsequences ${\cal S}_{{\cal E},(F_i)}$
of the system ${\cal S}_{\cal E}$, as described in Remark 2.A.3, in which 
the sequences $F_1\subset F_2\subset\dots$ of finite subtrees of $T$
satisfy the following conditions:
\item{(j1)} $F_1$ is a subtree of $T$ reduced to a single vertex;
\item{(j2)} for each $i\ge1$ the vertex set $V_{F_{i+1}}\setminus V_{F_i}$
consists of vertices adjacent to $F_i$ (i.e. of form $\omega(e)$ for some
$e\in N_{F_i}$).

\noindent
A choice of a sequence $(F_i)$ in the above way yields an inverse
sequence with the following properties:

\item{(p1)}
each manifold
$\hat K_{F_{i+1}}$ is obtained from $\hat K_{F_i}$ by means of  
operations of connected sum, with copies of the manifold $M$, 
performed at a finite set of pairwise
disjoint connecting disks in the manifold $\hat K_{F_i}$;

\item{(p2)} each map $h_{F_{i+1},F_i}$ maps each new copy of $M$
in $\hat K_{F_{i+1}}$, attached to $\hat K_{F_i}$ by means of
conneted sum at appropriate connecting disk, to this disk,
and it is the identity in the complement of the union
of the connecting disks.

\medskip\noindent
{\bf 2.C.8 Remark.} 
There are at least two examples in the literature where  
ideal boundaries of groups have been identified as certain trees of
manifolds, and this has been achieved
by referring to the description of these spaces as limits of
Jakobsche's inverse sequences, as described in Example 2.C.7. 

First, there is a paper [Fi] by Hanspeter Fischer, in which
the author identifies the CAT(0) boundaries of right angled
Coxeter groups having manifold nerves as trees of the corresponding
manifolds. A minor mistake in the statement of the main result
of [Fi] is corrected in Theorem 3.A.3 below.
The work of Fisher is complemented by the paper [P\'S],
by Piotr Przytycki and the author, where a vast class of trees
of manifolds, in dimensions $\le3$, is realized as Gromov
boundaries of Coxeter groups which are hyperbolic.

A different class of spaces, and associated hyperbolic groups,
is studied by Pawe{\l}  Zawi\'slak in [Z1]. He shows, 
among others, that Gromov boundary of a 7-systolic
3-dimensio\-nal orientable simplicial pseudomanifold
is the tree of 2-tori (known also as the Pontriagin sphere).
A related class of spaces is exhibited, for which the Gromov boundary
is the tree of projective planes.

\bigskip\noindent
{\bf 2.D More on the dimension of the limit of a tree system.}

\medskip
For the reasons explained in Subsection 1.D.3, we restrict
our attention to tree systems which have a universal finite upper bound
for the dimensions of the constituent spaces $K_t$.

\medskip\noindent
{\bf 2.D.1 Proposition.}
{\it Let $\Theta=(T,\{ K_t \},\{\Sigma_e \},\{ \phi_e \})$
be a tree system of metric compacta such that 
$$
\sup\{\dim(K_t):t\in V_T  \}=n<\infty.
$$
Suppose that for some constant $0\le c<1$ for each $e\in O_T$
there is a retraction $r_e:K_{\alpha(e)}\to\Sigma_e$
such that for each $e'\in N_{\alpha(e)}\setminus\{e\}$
the restriction $r_e|_{\Sigma_{e'}}$ is a $c$-contraction.
Suppose also that if the constant $c$ in the previous assumption is positive,
there exists another constant $C>0$ such that 
$\hbox{diam}(K_t)\le C$ for each $t\in V_T$.
Then $\dim(\lim\Theta)=n$.}

\medskip
Note that any family of maps $\{ r_e:e\in O_T \}$ as in 
Proposition 2.D.1, together with the trivial family
of extended spaces for $\Theta$ (see Example 2.A.1(1)), form
a contracting system ${\cal E}=(\{ \hat K_t=K_t \},\{ r_e \})$ 
of extended spaces and maps
for $\Theta$. Moreover, by Fact 2.B.3, $\cal E$ is then fine.
Thus, Proposition 2.D.1 is a special case of the following 
slightly more general result.
 
\medskip\noindent
{\bf 2.D.2 Proposition.}
{\it Let $\Theta=(T,\{ K_t \},\{\Sigma_e \},\{ \phi_e \})$
be a tree system of metric compacta such that 
$$
\sup\{\dim(K_t):t\in V_T  \}=n<\infty.
$$
Suppose that $\Theta$ admits a fine family $\cal E$
of extended spaces and maps, in which the corresponding 
family of extended spaces is trivial.
Then $\dim(\lim\Theta)=n$.}

\medskip\noindent
{\bf Proof of Proposition 2.D.2:}
Since we have an easy estimate
$$
\dim(\lim\Theta)\ge\sup\{\dim(K_t):t\in V_T  \}=n,
$$
it is sufficient to show the converse inequality $\dim(\lim\Theta)\le n$.

Let ${\cal S}_{\cal E}=(\{ \hat K_F \},\{ h_{FF'} \})$ be the inverse system
associated to $\cal E$. Since the family $\cal E$ is fine,
Theorem 2.B.4 implies that $\lim\Theta\cong \lim_\leftarrow {\cal S}_{\cal E}$. 
Since the family of extended spaces in $\cal E$
is trivial, for each finite subtree $F\subset T$ 
we have $\hat K_{F}=\#\{ K_t:t\in V_{F} \}$
(finite partial union), and hence
$\dim(\hat K_{F})=\max\{ \dim(K_t):t\in V_{F} \}\le n$. 
Moreover, by the properties of inverse limits we have
$\dim(\lim_{\leftarrow}{\cal S}_{\cal E})\le\sup_{F\in{\cal F}_T}\dim(\hat K_F)$.
This gives the required converse inequality
$\dim(\lim\Theta)\le n$, thus proving the proposition.

\bigskip\noindent
{\bf 2.D.3 Example: tree of internally punctured manifolds with boundary.}

We show how to apply Proposition 2.D.1 to calculate dimension of the limit
for the following class of tree systems
$\Theta=(T,\{K_t\}, \{ \Sigma_e \},\{ \phi_e\})$.
Fix any $n\ge1$, and suppose that for each $t\in V_T$
$$
K_t=M_t\setminus\bigcup\{ \hbox{int}(D):D\in{\cal D} \}
\hbox{\quad and \quad}
\{ \Sigma_e:e\in N_t \}=\{ \partial D:D\in{\cal D} \},
$$
where $M_t$ is a compact $n$-dimensional 
topological manifold with nonempty 
boundary, and $\cal D$ is a null and dense in $M_t$ family of
pairwise disjoint collared $n$-disks $D\subset\hbox{int}(M_t)$.
We call any $\Theta$ as above a
{\it dense tree system of internally punctured manifolds with boundary}.
We claim that for any such $\Theta$ the topological dimension
$\dim(\lim\Theta)=n-1$.

To prove the claim, fix any $t\in V_T$ and any $e\in N_t$, and 
suppose that $\Sigma_e=\partial D_0$, where $D_0\in{\cal D}$.
Observe that by collapsing all peripheral subsets
$\Sigma_{e'}:e'\in N_t\setminus\{ e \}$ to points one gets
the quotient space $K_t/\sim_e$ homeomorphic to 
$M_t\setminus\hbox{int}(D_0)$, via a homeomorphism which sends
$\partial M_t\cup\partial D_0\subset K_t/\sim_e$ identically to
$\partial M_t\cup\partial D_0\subset M_t\setminus\hbox{int}(D_0)$
(see Lemma 1.D.2.2).
We denote a homeomorphism 
$K_t/\sim_e\to M_t\setminus\hbox{int}(D_0)$ as above by $h_e$.
Observe also that, since $\partial M_t\ne\emptyset$, there exists
a retraction $g_e:M_t\setminus\hbox{int}(D_0)\to\partial D_0=\Sigma_e$.
(Existence of such a retraction is pretty obvious when $M_t$ is either smooth or PL, and it requires some effort in topological category; the latter case is carefully dealt with in [Z2].)
Thus, putting $r_e:=g_e\circ h_e$, we obtain a family $\{r_e\}$
of retractions as in the assumption of Proposition 2.D.1,
with constant $c=0$. Since we clearly have $\dim(K_t)=n-1$
for each $t$, Proposition 2.D.1 directly implies that
$\dim(\lim\Theta)=n-1$.

\medskip\noindent
{\bf 2.D.3.1 Remark.}
Let $\cal M$ be a dense tree system of internally punctured manifolds with
boundary, as defined above. Let $M$ be a connected topological manifold with nonempty
boundary, oriented or non-orientable, and suppose that all manifolds
$M_t$ used to describe $\cal M$ are homeomorphic to $M$.
Suppose also that, in case when $M$ is oriented, all connecting maps
$\phi_e$ of $\cal M$ respect orientations, i.e. they reverse the induced
orientations on the corresponding spheres. By the arguments as before
(see Examples 1.E.2 and 1.E.4) one easily shows that the system
$\cal M$ depends uniquely up to isomorphism on $M$ only.
We will call such system $\cal M$ the 
{\it dense tree system of internally punctured manifolds $M$},
and its limit $\lim{\cal M}$ the {\it tree of internally punctured
manifolds $M$}, denoted ${\cal X}_{int}(M)$.

\medskip\noindent
{\bf 2.D.4 Remark.} 
The following example shows that the equality as in the assertion of
Proposition 2.D.1 or 2.D.2 does not hold universally.
Let $\cal M$ be a dense tree system of closed oriented $(n+1)$-dimensional
manifolds $\{M_t\}$. Then the corresponding constituent spaces $K_t$
(obtained from the manifolds $M_t$ as
in Definition 1.E.2.2)
can be easily shown to have the topological dimension $\dim(K_t)=n$.
In particular, we have $\sup\{\dim(K_t):t\in V_T  \}=n$.
On the other hand, it is known that $\dim(\lim{\cal M})=n+1$, see 
Proposition (2.2) in [J2]. In the special case when all $M_t$
are $(n+1)$-spheres, this follows also from Lemma 1.D.2.1.

\bigskip
The next result concerns peripherally ANR tree systems $\Theta$, 
as defined in Section 2.C. It exhibits another class of examples for which
we have the equality 
$\dim(\lim\Theta)=\sup\{ \dim(K_t):t\in V_T \}$.

\medskip\noindent
{\bf 2.D.5 Proposition.}
{\it For any peripherally ANR tree system of metric compacta
$\Theta=(T,\{ K_t \},\{ \Sigma_e \},\{ \phi_e \})$ we have
$$
\dim(\lim\Theta)\le\max
\left(  \sup_{t\in V_t}\dim(K_t),\sup_{e\in O_T}\dim(\Sigma_e)+1 \right).
$$
In particular, if $\Theta$ is peripherally ANR,  
if $\sup_{t}\dim(K_t)=n<\infty$ and $\sup_e\dim(\Sigma_e)\le n-1$,
then $\dim(\lim\Theta)=n$.}

\medskip\noindent
{\bf Proof:}
Consider any conical and 0-contracting family $\cal E$ of extended spaces
and maps for $\Theta$, the existence of which is justified by Proposition 2.C.2,
and recall that it is fine. Let ${\cal S}_{\cal E}$ be the associated inverse system
for $\Theta$ induced by $\cal E$. 
For each $t\in V_T$ the family $\{ K_t \}\cup\{ \Delta_e:e\in N_t \}$
is a countable closed covering of $\hat K_t$ and hence,
by the countable union theorem (see [Eng, Theorem 7.2.1]),
we have 
$\dim(\hat K_t)=\max\left(\dim(K_t),\sup_{e\in N_t}\dim(\Delta_e)\right)$.
Since the family $\cal E$ is conical,
and since we have the equality $\dim(cone(X))=\dim(X)+1$ for any
compact metric space, it follows that
$$
\dim(\hat K_t)=\max\left(\dim(K_t),\sup_{e\in N_t}\dim(\Sigma_e)+1\right).
$$
Using this, and applying once again the countable union theorem
(this time to a finite union),
we get the following estimate for each finite subtree $F$ of $T$:
$$
\dim(\hat K_F)\le
\left(  \sup_{t\in V_t}\dim(K_t),\sup_{e\in O_T}\dim(\Sigma_e)+1 \right).
$$
Finally, since by Theorem 2.B.4 we have 
$\lim\Theta\cong\lim_\leftarrow {\cal S}_{\cal E}$,
and since by the properties of inverse limits we have
$\dim(\lim_{\leftarrow}{\cal S}_{\cal E})\le\sup_{F\in{\cal F}_T}\dim(\hat K_F)$,
we get the required estimate for $\dim(\lim\Theta)$, as in the first
assertion of the proposition.

The second assertion follows from the first one,
and from the inequality (1.D.3).

\bigskip\bigskip
\centerline{\bf 3. Modifications of tree systems.}

\medskip
In this part of the paper we describe some natural operations
on tree systems which do not affect their limits. We also present
several applications of these operations for justifying various
properties of trees of manifolds. We are convinced that these operations
provide a powerful tool for the future study of more general classes
of trees of spaces.

\bigskip\noindent
{\bf 3.A Consolidation of a tree system.}

\medskip
We describe an operation which turns one tree system of spaces
into another by merging the constituent spaces of the initial system, 
and forming a new system out of bigger pieces (corresponding to
a family of pairwise disjoint subtrees in the underlying tree
of the initial system).
As we show below (see Theorem 3.A.1), 
this operation does not affect the limit
of a system.

Let $\Theta=(T,\{ K_t \},\{ \Sigma_e \},\{ \phi_e \})$ be a tree
system of metric compacta.
Let $\Pi$ be a {\it partition} of a tree $T$ {\it into subtrees},
i.e. a family of subtrees $S\subset T$ such that the vertex sets
$V_S:S\in{\Pi}$ are pairwise disjoint and cover all of $V_T$.
We allow that some of the subtrees
$S\in{\Pi}$ are just single vertices of $T$.

We define a {\it consolidation of $\Theta$ with respect to $\Pi$},
denoted $\Theta_{\Pi}$, to be the following tree system 
$(T_{\Pi}, \{ K_S \}, \{ \Sigma_e \}, \{ \phi_e \})$.
As a vertex set $V_{{\Pi}}$ of the tree $T_{{\Pi}}$ we take the family $\Pi$,
and as the edge set $O_{{\Pi}}$ the set 
$\{ e\in O_T:e\notin\cup_{S\in{\Pi}}\,O_S \}$. 
Clearly, for each oriented edge $e\in O_{\Pi}$ the initial
vertex $\alpha_{\Pi}(e)$ is this subtree $S\in{\Pi}$
for which $\alpha(e)\in V_S$.
Similarly, $\omega_{\Pi}(e)$ is this subtree $S\in{\Pi}$
for which $\omega(e)\in V_S$.

For any subtree $S\in{\Pi}$ denote by $\Theta_S$ the restriction
of $\Theta$ to $S$, and put $K_S:=\lim\Theta_S$. Note that for any
$e\in O_{\Pi}$ we have the canonical inclusion 
$\Sigma_e\subset K_{\alpha_{\Pi}(e)}$, and if we put
$N_S=\{ e\in O_{\Pi}:\alpha_{\Pi}(e)=S \}$, the family
$\Sigma_e:e\in N_S$ of subsets of $K_S$ consists of pairwise
disjoint sets. Moreover, by Proposition 1.C.1, the subsets
$\Sigma_e\subset K_{\alpha_{\Pi}(e)}$ are all compact
and each family $\Sigma_e:e\in N_S$ is null.
This justifies that the just described tuple 
$$
\Theta_{\Pi}:=(T_{\Pi}, \{ K_S:S\in{\Pi} \}, \{ \Sigma_e:e\in O_{\Pi} \},
\{ \phi_e:e\in O_{\Pi} \})
$$
is a tree system of metric compacta.

\medskip\noindent
{\bf 3.A.1 Theorem.}
{\it For any tree system $\Theta$ of metric compacta,
and its any consolidation $\Theta_{\Pi}$, the limits
$\lim\Theta$ and $\lim\Theta_{\Pi}$ are canonically
homeomorphic.}

\bigskip\noindent
{\bf Proof:} 
We first describe the natural map $i_\Pi:\lim\Theta_\Pi\to\lim\Theta$,
and then show it is a homeomorphism. 

To define $i_\Pi(x)$ for all $x\in\lim\Theta_\Pi$,
we consider three possible positions of $x$ in $\lim\Theta_\Pi$.
First, suppose that $x\in K_t\subset K_S\subset\lim\Theta_\Pi$,
for some $t\in V_S$ and some $S\in\Pi$.
View $x$ as an element of $\lim\Theta$ via the canonical inclusion
$K_t\subset\lim\Theta$, and put $i_\Pi(x)=x$.
Second, suppose that $x\in E_S\subset\lim\Theta_S= K_S\subset\lim\Theta_\Pi$
is an end of some tree $S\in\Pi$. Since we have the canonical inclusion
$E_S\subset E_T$, $x$ is also an end of $T$, and hence an element
of $\lim\Theta$. We put again $i_\Pi(x)=x$. Finally, in the remaining case
$x$ is an end of the tree ${T_\Pi}$. Suppose that $x$ is represented
by a ray $\rho=(S_0,S_1,\dots)$ in the tree $T_\Pi$. Then $\rho$
determines the ray $\rho_T$ in $T$ by
$$
\rho_T=(e_1,[\omega(e_1),\alpha(e_2)],
e_2,[\omega(e_2),\alpha(e_3)],e_3,\dots),
$$
where $e_i$ is the edge in $T$ connecting $S_{i-1}$ to $S_i$,
and $[\omega(e_i),\alpha(e_{i+1})]$ are the paths in $T$
(sometimes perhaps empty) connecting the corresponding vertices.
Put $i_\Pi(x)=[\rho_T]\in E_T\subset\lim\Theta$.
We skip the straightforward verification that
$i_\Pi$ is well defined, and that it is a bijection.

Since any continuous bijection between compact metric spaces 
is a homeomorphism,
to prove that $i_\Pi$ is a homeomorphism, it is sufficient to show that
it is continuous. To do this, we will show that for any open set
$V\subset\lim\Theta$ from the basis $\cal B$ the preimage
$i_\Pi^{-1}(V)$ is open in $\lim\Theta_\Pi$.
For this we need two claims.

\medskip\noindent
{\bf Claim 1.} 
{\it For any $S\in\Pi$ the restriction $i_\Pi|_{K_S}$ is a homeomorphism
on its image (a topological embedding).}

\medskip
Since $K_S$ is compact and $i_\Pi$ is injective, 
to prove Claim 1 it is sufficient to show that 
$i_\Pi|_{K_S}$ is continuous. Let $G(U)\in{\cal B}$, where
$U$ is an open subset in $K_F$ for some finite $F\subset T$,
and $U$ is ${\cal A}_F$-saturated. We need to show that
$K_S\cap i_\Pi^{-1}(G(U))$ is open in $K_S$.
Suppose first that $F\cap S=\emptyset$. It is not hard to realize
that in this case $K_S\cap i_\Pi^{-1}(G(U))$
is either empty or coincides with $K_S$,
and thus the assertion follows.
Suppose then that $F\cap S\ne\emptyset$. Viewing $K_{F\cap S}$
as subset in $K_F$, put $U_S:=U\cap K_{F\cap S}$,
and note that $U_S$ is open in $K_{F\cap S}$ and saturated
with respect to the family ${\cal A}_{F\cap S}^{\Theta_S}$ of peripheral subsets
in $K_{F\cap S}$ viewed as a partial union of the system
$\Theta_S$. Moreover, it is not hard
to observe that $K_S\cap i_\Pi^{-1}(G(U))=G(U_S)\in{\cal B}_S$,
where ${\cal B}_S$ is the standard basis in the limit $K_S$
of the system $\Theta_S$ (as described in Section 1.C).
Thus Claim 1 follows.

\medskip\noindent
{\bf Claim 2.} 
{\it For any finite subtree $F_0\subset T_\Pi$ the restriction 
$i_\Pi|_{K_{F_0}}$ is a homeomorphism on its image.}

\medskip
Since $K_{F_0}$ is a finite union of its compact subsets of form $K_S$,
it is compact itself. Moreover, in view of Claim 1, the same fact implies
that $i_\Pi|_{K_{F_0}}$ is continuous. Since this map is also injective,
the assertion of Claim 2 follows.

\medskip
Coming back to the proof that $i_\Pi$ is continuous,
let $V=G(U)$ for some open and ${\cal A}_F$-saturated subset $U\subset K_F$, 
where $F\subset T$ is some finite subtree.
Let $F_\Pi$ be the subtree of $T_\Pi$ spanned by the
vertices represented by those subtrees $S\in\Pi$ which
intersect $F$; clearly, $F_\Pi$ is finite. 
Note that, by Claim 2, the set
$U_\Pi:=K_{F_\Pi}\cap i_\Pi^{-1}(V)$ is open in $K_{F_{\Pi}}$
(because $i_\Pi(K_{F_\Pi})\cap V$ is open in $i_\Pi(K_{F_\Pi})$).
Moreover, since $U$ is ${\cal A}_F$-saturated, it is not hard
to realize that $U_\Pi$ is ${\cal A}_{F_\Pi}$-saturated, where
${\cal A}_{F_\Pi}$ is the appropriate family of peripheral
subsets of the system $\Theta_\Pi$ in its partial union $K_{F_\Pi}$.
Finally, observe that $i_\Pi^{-1}(V)=G(U_\Pi)\in{\cal B}_\Pi$
(${\cal B}_\Pi$ denotes the standard basis for the topology in $\lim\Theta_\Pi$),
and hence this preimage is open in $\lim\Theta_\Pi$.
This finishes the proof.

\medskip
The next example illustrates how one can apply the procedure
of consolidation to justify that limits of various classes of tree
systems are homeomorphic.

\medskip\noindent
{\bf 3.A.2 Example: dense trees of finite families of manifolds.}

Let ${\cal N}=\{ M_1,\dots,M_k \}$ be a finite family of closed
connected oriented topological manifolds of the same dimension. Let
${\cal M}=(T,\{ K_t\},\{ \Sigma_e \},\{ \phi_e \})$ be a dense
tree system of manifolds from $\cal N$. For each $t\in V_T$
let $i_t\in\{ 1,\dots,k \}$ be this index for which the space $K_t$
has a form $M_{i_t}^\circ$, as described in Section 1.E, just
after Definition 1.E.2.2. We say that $\cal M$ is {\it 2-saturated}
if for each $t\in V_T$ and every $j\in\{ 1,\dots,k \}$
there are at least two distinct edges $e\in N_t$ such that
$i_{\omega(e)}=j$.

If $\cal M$ is 2-saturated, it is not hard
to construct a partition $\Pi$ of $T$ such that for each $S\in\Pi$
the vertex set $V_S$ contains exactly one vertex $s$ with $i_s=j$,
for each $j\in\{1,\dots,k \}$ (in particular, the cardinality of
each $V_S$ is $k$). Consider the consolidated tree system 
${\cal M}_\Pi$ for the partition $\Pi$. It is easy to note that the
constituent spaces of ${\cal M}_\Pi$ are all of form
$(M_1\#\dots\#M_k)^\circ$, and that ${\cal M}_\Pi$ is
(isomorphic to) the dense tree system of manifolds $M_1\#\dots\#M_k$.
In view of Theorem 3.A.1 and Proposition 1.E.2.3 we get the following.

\medskip\noindent
{\bf 3.A.2.1 Proposition.}
{\it For any 2-saturated dense tree
system $\cal M$ of manifolds from $\cal N$
the limit $\lim{\cal M}$ is homeomorphic to the Jakobsche space
${\cal X}(M_1\#\dots \#M_k)$.}
 
\medskip
The above observation can be generalized as follows. Let $\cal N$
be a family of manifolds as above, and let $\mu=(m_1,\dots,m_k)$
be a tuple of arbitrary positive integers. If $\cal M$ is a tree system
as above, one can easily construct a partition $\Pi_\mu$ 
such that for each $S\in\Pi_\mu$ and each $j\in\{ 1,\dots,k \}$
the vertex set $V_S$ contains exactly $m_j$ elements $s$ with
$i_s=j$. By the argument as above, we get that 
$\lim{\cal M}=\lim{\cal M}_{\Pi_\mu}$
is homeomorphic to the Jakobsche space
${\cal X}(m_1M_1\#\dots\#m_kM_k)$, where
each $m_iM_i$ is the connected sum of $m_i$ copies of $M_i$.
As a consequence, we get the following.

\medskip\noindent
{\bf 3.A.2.2 Corollary.}
{\it Let $M_1,\dots,M_k$ be any closed connected oriented topological manifolds
of the same dimension, and let $m_1,\dots,m_k$ be arbitrary
positive integers. Then the Jakobsche spaces
${\cal X}(M_1\#\dots\#M_k)$ and
${\cal X}(m_1M_1\#\dots\#m_kM_k)$ are homeomorphic.
}

\medskip
Now, let ${\cal K}=\{ N_1,\dots,N_k \}$ be a family of closed
connected topological
manifolds of the same dimension, at least one of which
is non-orientable. We also assume that the orientable manifolds
in $\cal K$ have no distinguished orientation. By a {\it tree system
of manifolds from $\cal K$} we mean a tree system of manifolds
in which every constiuent space has a form $N^\circ$ for some
$N\in{\cal K}$, and in which the connecting maps are arbitrary
homeomorphisms between the corresponding spherical peripheral
subspaces. Note that a connected sum $N_1\#\dots\#N_k$
(again, with arbitrary connecting homeomorphisms) is unique
up to homeomorphism (this is a consequence of Lemma 1.E.4.1).

The arguments as above, enriched by application of Lemma 1.E.4.1,
yield the following variation on Proposition 3.A.2.1 and Corollary
3.A.2.2.

\medskip\noindent
{\bf 3.A.2.3 Proposition.}
{\it Let ${\cal K}=\{ N_1,\dots,N_k \}$ be a family of closed
connected topological manifolds of the same dimension, at least one of which
is non-orientable. Let $\cal M$ be any 2-saturated dense tree system 
of manifolds from $\cal K$. Then the limit $\lim{\cal M}$
is homeomorphic to the Jakobsche space 
${\cal X}(N_1\#\dots\# N_k)$ (as defined in Example 1.E.4).
Moreover, for any positive integers $m_1,\dots,m_k$ the space
${\cal X}(m_1N_1\#\dots\# m_kN_k)$ is also homeomorphic to 
${\cal X}(N_1\#\dots\# N_k)$.}

\medskip\noindent
{\bf 3.A.2.4 Remark.}
In [J2] Jakobsche considered inverse sequences related to
dense tree systems of manifolds from a family $\cal N$, and established
a much weaker variant of Proposition 3.A.2.1. Namely, using our terms,
he got uniqueness up to homeomorphism under the following
(much stronger than 2-saturation) assumption for a dense tree system
of manifolds from $\cal N$: {\it for each $M\in {\cal N}$ and for
each $t\in T$ the set of collared disks 
$\Delta_e\subset M_t:M_{\omega(e)}\cong M$ is dense in $M_t$}
(see Theorem (4.6) in [J2]).
His proof consisted of showing something equivalent to the fact
that any two tree systems satisfying the above assumption are isomorphic.
Observe that,
under our weaker assumption of 2-saturation, tree systems as in
Proposition 3.A.2.1 needn't be isomorphic, so the proof of this
proposition necessarily requires some argument involving modifications
of tree systems.


\medskip
As another application of the technique of consolidation we present
the
following correction to the main result from the paper [Fi]
by Hanspeter Fischer.

\medskip\noindent
{\bf 3.A.3 Theorem.}
{\it Suppose $W$ is a right angled Coxeter group whose nerve
is a flag PL triangulation of a closed oriented manifold $M$,
and let $\overline M$ be the same manifold with reversed orientation.
Then the CAT(0) boundary of $W$ (i.e. the boundary of the
Coxeter-Davis complex for $W$) is the Jakobsche space
${\cal X}(M\#\overline M)$.}

\medskip
In the original (wrong) statement of this result in [Fi] 
instead of the space
${\cal X}(M\#\overline M)$, there appears
the space ${\cal X}(M)$, which is in general different
(for example, it is not hard to show that the spaces
${\cal X}(CP^2)$ and ${\cal X}(CP^2\#\overline{CP}^2)$
are not homeomorphic, by referring to the properties
of their \v Cech cohomology rings).

\medskip\noindent
{\bf Proof of Theorem 3.A.3:} 
We indicate a necessary minor modification of the 
argument provided in [Fi].
The author argues by showing that $\partial W$ 
is homeomorphic to the inverse limit of some inverse sequence
of manifolds of the form as in Example 2.C.7.
A part of his argument which requires correction
is this (see the beginning of the proof of Theorem 3.7 in [Fi]).
Let $X_W$ be the Coxeter-Davis complex for $W$, obtained as the 
union $X_W=\cup\{ gQ:g\in W \}$, where $Q$ is the Davis' cell
for $W$ (topologically equal to the cone over the manifold $M$).
Let $|g|$ be the word norm for elements $g\in W$ with respect
to the standard generating set. Let $X_k=\cup\{ gQ:|g|\le k \}$
and let $M_k=\partial X_k$. The author claims that $M_k$ is the
connected sum of the appropriate number of copies of $M$,
but it is clear from the way $X_k$ is formed out of copies of $Q$
that in fact $M_k$ is the connected sum of copies of both $M$
and $\overline M$. Moreover, since any two adjacent copies of $Q$
in $X_W$ are obtained from one another by reflection,
$M_{k+1}$ is obtained from $M_k$ by connected sum with
copies of $M$ if $k$ is odd, and with copies of $\overline M$
if $k$ is even (in particular, $M_0=\partial Q=M$).
From this, using the remaining arguments of Fischer,
one deduces that $\partial W$ is homeomorphic to the limit
of a 2-saturated dense tree system of manifolds $M$ and $\overline M$.
By Proposition 3.A.2.1, this yields the assertion.

\bigskip\noindent
{\bf 3.B Tree decomposition of a compact metric space.}

\medskip
In Sections 3.B and 3.C we describe an operation on a tree system
that is inverse to that of consolidation (as described in Section 3.A).
In this section, we start with an elementary case when the initial
tree system is trivial, i.e. the underlying tree $T$ is reduced to
a single vertex. The operation is then called 
{\it tree decomposition} of a compact metric space.

We start with introducing terminology related to the concept
of {tree decomposition}.
The main result of the section, which relates this concept
to that of a tree system and its limit, is Theorem 3.B.10 below.

At the end of the section we show how to use tree decompositions
to prove that limits of some tree systems are homeomorphic to some
explicit spaces (Example 3.B.12).

\medskip\noindent
{\bf 3.B.1 Definition.} 
An {\it elementary splitting} of a compact metric space $K$
is a triple $(A,\{Y,Z\})$ of compact subspaces of $K$ such that
$Y\cup Z=K$, $Y\cap Z=A$, and $A$ is a nonempty 
proper subset both in $Y$ and in $Z$. 
The set $A$ is called the {\it separator}
of the splitting, and the sets $Y,Z$ are the {\it halfspaces}. 
Moreover, the sets $Y\setminus A$ and $Z\setminus A$ will be called
the {\it open halfspaces of the splitting}. If $H$ is any halfspace of the
splitting above (i.e. $H=Y$ or $H=Z$), we denote by $\dot{H}$
the corresponding open halfspace, and by $H^c$ the {\it complementary}
(or {\it opposite}) halfspace (equal to $K\setminus\dot{H}$).

\medskip
Note that for any splitting as above we have the following:
\item{(1)} the set $K\setminus A$ is disconnected and the open
halfspaces are the unions of connected components in this set;
\item{(2)} $K$ is canonically homeomorphic to the space
$Y\cup_A Z$ obtained from the disjoint union of $Y$ and $Z$
by gluing through the identity on $A$; equivalently, $K$ is the limit
of the tree system whose underlying tree $T$ is a single edge, 
the constituent spaces $K_t$ are $Y$ and $Z$, the peripheral subspaces $\Sigma_e$ coincide with $A$,
and the connecting maps $\phi_e$ are the identities on $A$.

\medskip\noindent
{\bf 3.B.2 Definition.}
Given two elementary splittings $(A_i,\{ Y_i,Z_i \})$, $i=1,2$, of
a compact metric space $K$, we say they {\it do not cross} if
for at least one pair of halfspaces $H_1,H_2$ 
selected from those splittings we have $H_1\cap H_2=\emptyset$.

\medskip\noindent
{\bf 3.B.3 Remark.}
Note that the noncrossing condition $H_1\cap H_2=\emptyset$
has the following consequences:
{(1)} $A_1\cap A_2=\emptyset$,
{(2)} $H_1\subset \dot{H}_2^c$ and $H_2\subset \dot{H}_1^c$.

\medskip\noindent
{\bf 3.B.4 Definition.} Given three pairwise noncrossing splittings
$(A_i,\{ Y_i,Z_i \})$, $i=1,2,3$, of
a compact metric space $K$, we say that $A_2$ separates $A_1$ from $A_3$
if for appropriately chosen halfspace $H$ for $A_2$ we have
$A_1\subset\dot{H}$ and $A_3\subset\dot{H}^c$.

\medskip
We will be interested in countable (usually infinite) families of
pairwise noncrossing splittings satisfying some additional properties,
which we now describe.

\medskip\noindent
{\bf 3.B.5 Definition.}
Let ${\cal C}=(A_\lambda,\{ Y_\lambda,Z_\lambda \})_{\lambda\in\Lambda}$ 
be a family
of pairwise noncrossing splittings of a compact metric space $K$. 
We say that $\cal C$ is {\it discrete} if for any two separators
$A,A'$ from $\cal C$ there is only finitely many separators in $\cal C$
that separate $A$ from $A'$. 

\medskip
Any discrete family $\cal C$ of pairwise noncrossing splittings of $K$
determines the associated family of {\it domains obtained by splitting},
and the {\it dual tree} $T_{\cal C}$. 
We start with describing the domains.

Consider any separator $A$ from $\cal C$ and any halfspace $H$
related to $A$. The pair $(A,H)$ determines a domain $\Omega_{A,H}$ 
described as follows.
Let $A=A_{\lambda_0}$, and put 
$\Lambda_0=\Lambda\setminus\{ \lambda_0 \}$.
For any separator $A_\lambda$ with $\lambda\in\Lambda_0$ 
consider this halfspace $H_\lambda$
which contains $A$. We have to consider the following two cases. 
First, suppose that for each 
$\lambda\in\Lambda_0$ we have $H\subset H_\lambda$.
We then put $\Omega_{A,H}:=H$ and note that this domain is disjoint with
all separators from $\cal C$ other than $A$.
In the second case, there is some $\lambda\in\Lambda_0$ with 
$H\subset\hskip-9pt/ \,\,H_\lambda$. By discreteness, there is also
such $\lambda$ for which no separator from $\cal C$ separates
$A$ from $A_\lambda$. Denote by $\Lambda_{A,H}$ the set of all
$\lambda\in\Lambda_0$ for which $H\subset\hskip-9pt/ \,\,H_\lambda$
and no separator from $\cal C$ separates
$A$ from $A_\lambda$. By what was said above, this set is nonempty.
Put 
$$
\Omega_{A,H}:=H\cap\bigcap_{\lambda\in\Lambda_{A,H}}H_\lambda
$$
and note that this set satisfies the following properties:
\item{(1)} $\Omega_{A,H}$ is compact, contains the sets $A$
and $A_\lambda:\lambda\in\Lambda_{A,H}$, and it is disjoint
with all other separators from $\cal C$;
\item{(2)} for each $\lambda\in\Lambda_{A,H}$
we have $\Omega_{A_\lambda,H_\lambda}=\Omega_{A,H}$.

\noindent
A {\it domain} in $K$ induced by $\cal C$ is any subset $\Omega=\Omega_{A,H}$ as above. A domain $\Omega$ is called {\it adjacent} to a separator $A$ of $\cal C$ if it contains $A$ or equivalently, if $\Omega=\Omega_{A,H}$ for some halfspace $H$ related to $A$. For each separator $A$ of $\cal C$ there are exectly two domains induced by $\cal C$ and adjacent to $A$. 

We are now ready to describe the dual tree $T_{\cal C}$ of a discrete
family of pairwise noncrossing splittings $\cal C$.
As a vertex set $V_{\cal C}$ 
we take the set of all domains for $\cal C$, as described
above. As the set $O_{\cal C}$ of orieted edges we take
the set of all pairs $(A,H)$ as above, and we put
$\omega_{\cal C}(A,H)=\Omega_{A,H}$ and $\overline{(A,H)}=(A,H^c)$. 
It is an easy exercise to show
that the graph obtained in the above way is connected, and contains
no loops, thus being a tree. We denote this tree by $T_{\cal C}$
and call the {\it dual tree of $\cal C$}.

\medskip\noindent
{\bf Remark.}
Note that for any discrete family $\cal C$ of pairwise noncrossing
splittings  of a compact metric space $K$ the corresponding
dual tree $T_{\cal C}$ is locally countable (and hence countable).
More precisely, given a vertex $\Omega$ of $T_{\cal C}$ (i.e. a domain
for $\cal C$), consider the set $N_\Omega$ of all edges $(A,H)$ 
in $T_{\cal C}$ satisfying $\alpha_{\cal C}(A,H)=\Omega$. Note that the
corresponding family of open halfspaces in $K$, 
$\{ \dot H :(A,H)\in N_\Omega\}$ consists of pairwise disjoint
nonempty open subsets of $K$. Since any compact metric space is 
separable, it follows that the set $N_\Omega$ is countable,
thus justifying the remark.

\medskip
The above construction of the dual tree motivates the first part 
of the following.

\medskip\noindent
{\bf 3.B.6 Definition.}
\item{(1)} A discrete family $\cal C$ of pairwise
noncrossing splittings of a compact metric space $K$ will be called
{\it a tree decomposition} of $K$.
\item{(2)} A tree decomposition 
${\cal C}=
\{ (A_\lambda,\{ Y_\lambda,Z_\lambda \}) \}_{\lambda\in\Lambda}$
of $K$ is {\it fine} if for each $\epsilon>0$ the set
$
\{ \lambda\in\Lambda:
\min[\hbox{diam}(Y_\lambda),\hbox{diam}(Z_\lambda)]>\epsilon \}
$
is finite.

\medskip\noindent
{\bf 3.B.7 Example.}
Let $K=\lim\Theta$ be the limit of a tree system $\Theta$ of metric compacta,
and suppose that it is {\it essential}, in the sense that for any $e\in O_T$
the set $\Sigma_e$ is a proper subset in $K_{\alpha(e)}$.
For any edge $e\in O_T$, 
let $T_e$ denote the maximal subtree of $T\setminus\hbox{int}(|e|)$
that contains $\omega(e)$.
Put $H_e$ to be the limit of the restricted system $\Theta_{T_e}$.
Note that for each $e\in O_T$ the triple $(\Sigma_e,\{ H_e,H_{\bar e} \})$,
viewed as consisting of subsets of $K$, is an elementary splitting of $K$.
Note also that the family of all splittings of $K$ having this form
yields a tree decomposition of $K$. We denote this decomposition
by ${\cal C}(\Theta)$. Moreover, this decomposition is fine, which
can be deduced by the same methods as in the proof of the last assertions in
Proposition 1.C.1 (as given in Section 1.D).

\medskip
One of the consequences of fineness of a tree decomposition is the
following.

\medskip\noindent
{\bf 3.B.8 Fact.}
{\it Let $\cal C$ be a fine tree decomposition of a compact metric space $K$.
Then for any domain $\Omega$ for $\cal C$ the family ${\cal A}_\Omega$ 
of all separators from $\cal C$ adjacent to $\Omega$ is null.}

\medskip\noindent
{\bf Proof:} Note that ${\cal A}_\Omega$ coincides with the family
of those separators $A$ from $\cal C$ for which there is a halfspace $H$
(related to $A$) such that
$\alpha_{\cal C}(A,H)=\Omega$. Since any
opposite halfspace $H^c$ to a halfspace $H$ as in the previous 
sentence contains $\Omega$, and since for such $H^c$ we have
$\hbox{diam}(H^c)\ge\hbox{diam}(\Omega)>0$,
it follows from fineness of $\cal C$ that the family of halfspaces 
$\{ H:\alpha_{\cal C}(A,H)=\Omega \}$ 
is null. Consequently, since we have inclusions $A\subset H$,
the family ${\cal A}_\Omega$ is also null.

\medskip\noindent
{\bf 3.B.9 Definition.}
Given a fine tree decomposition $\cal C$ of a compact metric space $K$, 
{\it the tree system $\Theta_{\cal C}$ associated to $\cal C$} is described 
as follows. The underlying tree for 
$\Theta_{\cal C}$ is the dual tree $T_{\cal C}$.
For each vertex $t\in V_{\cal C}$ represented by some domain $\Omega$
we put $K_t=\Omega$. For each oriented edge $e=(A,H)\in O_{\cal C}$
we put $\Sigma_e=A$ and $\phi_e=id_A$. In view of Fact 3.B.8, 
this well defines a tree system of metric compacta (which is moreover
essential).

\medskip\noindent
{\bf Remark.}
It is not hard to realize that if ${\cal C}(\Theta)$ is the tree decomposition
of $\lim\Theta$ described in Example 3.B.7 then the associated tree system
$\Theta_{{\cal C}(\Theta)}$ is canonically isomorphic with $\Theta$.

\medskip
The main result of this section is the following.

\medskip\noindent
{\bf 3.B.10 Theorem.}
{\it For any fine tree decomposition $\cal C$ of a compact metric space $K$
the limit $\lim\Theta_{\cal C}$ of the associated tree system is canonically
homeomorphic to $K$.}

\medskip
To prove Theorem 3.B.10 we need the following result which exhibits
consequences of fineness.

\medskip\noindent
{\bf 3.B.11 Fact.}
{\it Let $\cal C$ be a fine tree decomposition of a compact metric space $K$.}
\item{(1)} 
{\it For each domain $\Omega$ for $\cal C$ the family of halfspaces
from $\cal C$ given by
${\cal H}_\Omega:=\{ H:\alpha_{\cal C}(A,H)=\Omega \}$ is null.}
\item{(2)} 
{\it For each ray $\rho=(e_1,e_2,\dots)$ in the dual tree $T_{\cal C}$,
with $e_i=(A_i,H_i)$, the corresponding family $\{ H_i \}$ of halfspaces
is null in $K$.}

\medskip\noindent
{\bf Proof:} The proof of part (1) has already been given
in the proof of Fact 3.B.8. To prove part (2),
put $\Omega=\alpha_{\cal C}(e_1)$ and note that
for each $i$ we have $\Omega\subset H_i^c$. Consequently,
we have $\hbox{diam}(H_i^c)\ge\hbox{diam}(\Omega)>0$.
It follows then from fineness of $\cal C$ that
$\lim_{i\to\infty}\hbox{diam}(H_i)=0$, which completes
the proof.

\medskip\noindent
{\bf Proof of Theorem 3.B.10:}
 We use the identification of $\lim\Theta_{\cal C}$
with the inverse limit $\lim_{\longleftarrow}{\cal S}_{\Theta_{\cal C}}$, 
provided by Proposition 1.D.1.

We describe
a natural map $\mu:K\to \lim\Theta_{\cal C}$.
Denoting 
$\Theta_{\cal C}=(T_{\cal C},\{ K_t \},\{ \Sigma_e \},\{ \phi_e \})$,
for each finite subtree $F$ of $T_{\cal C}$ consider a map 
$\mu_F:K\to K_F^*$ (where $K_F^*$ is the corresponding reduced partial
union of the system $\Theta_{\cal C}$) defined as follows. Denoting $e=(A_e,H_e)$
for all $e\in O_{T_{\cal C}}$, we have that $K_F=\bigcap_{e\in N_F}H_e^c$
and $K_F^*=K_F/\{ A_e:e\in N_F \}=K/\{ H_e:e\in N_F \}$.
The latter equality is a homeomorphism by nullness of the family
$H_e:e\in N_F$ (which follows from Fact 3.B.11(1) by observing
that $H_e:e\in N_F$ is a subfamily of the union of finitely many
families as in this fact).
Thus, we take as $\mu_F$ the corresponding quotient map
$K\to K/\{ H_e:e\in N_F \}=K_F^*$. Since the maps $\mu_F$
commute with the maps in the inverse system ${\cal S}_{\Theta_{\cal C}}$,
they induce a continuous map 
$\mu:K\to \lim_{\longleftarrow}{\cal S}_{\Theta_{\cal C}} =\lim\Theta_{\cal C}$.
Using the arguments as in the proof of Proposition 1.D.1,
we verify that $\mu$ is a bijection. (In the proof of injectivity
one needs to use Fact 3.B.11(2)).
Since both spaces $K$ and $\lim\Theta_{\cal C}$ are compact,
$\mu$ is a homeomorphism,
which completes the proof.

\medskip\noindent
{\bf 3.B.12 Example: boundary tree of disks is a disk.}

For any $n\ge2$, let $D^n$ be the $n$-disk, and let $\cal D$
be a dense and null family of pairwise disjoint collared $(n-1)$-disks
in the boundary sphere $\partial D^n$.
Note that, due to Toru\'nczyk's Lemma 1.E.2.1 (followed by Alexander's
trick), the tuple $(D^n,{\cal D})$ is unique up to homeomorphism.
Consider the unique tree system in which all constituent spaces
are homeomorphic to $D^n$, and all families of peripheral subspaces
coincide, up to ambient homeomorphism, with $\cal D$.
We denote this tree system by ${\cal M}_\partial(D^n)$ and call it the
{\it dense boundary tree system of $n$-disks}. 

\medskip\noindent
{\bf 3.B.12.1 Lemma.}
{\it The limit $\lim {\cal M}_\partial(D^n)$ is the $n$-disk.}

\medskip\noindent
{\bf Proof:}
We will describe a fine tree decomposition $\cal C$ of the $n$-disk
$D^n$ such that the associated tree system $\Theta_{\cal C}$
is isomorphic to ${\cal M}_\partial(D^n)$. By applying Theorem
3.B.10, and in view of Lemma 1.E.1, this gives the assertion.

To construct $\cal C$, view $D^n$ as the standard $n$-disk in $R^n$,
and consider the group $\hbox{\it M\"ob}(D^n)$ of all M\"obius 
transformations of $D^n$, i.e. those M\"obius transformations of $R^n$
which preserve $D^n$. Viewing $\hbox{int}(D^n)$ as the Poincare
disk model for the hyperbolic $n$-space, we will think of 
$\hbox{\it M\"ob}(D^n)$ as the group of all hyperbolic isometries
in its action on the completion of the hyperbolic $n$-space
by its ideal boundary.

Let $\cal D$ be any null and dense family of pairwise disjoint round
$(n-1)$-disks in $\partial D^n$.
For each $D\in{\cal D}$, let $H_D$ be the hyperbolic halfspace in $D^n$
such that $H_D\cap\partial D^n=D$, and let $H_D^c$ be the
opposite halfspace. Denote by $A_D$ the hyperplane
in $D^n$ bounding $H_D$, and by $s_D$ the heperbolic reflection
with respect to $A_D$, which clearly belongs to $\hbox{\it M\"ob}(D^n)$.

Let $\Gamma$ be the subgroup of $\hbox{\it M\"ob}(D^n)$ generated
by all elements $s_D:D\in{\cal D}$. Obviously, $\Gamma$ is then
an infinitely generated free reflection group with the fundamental domain
$$
\Omega_0:=\bigcap_{D\in{\cal D}}H_D^c.
$$ 
Algebraically, $\Gamma$ is the free product of its order 2 subgroups
$\langle s_D \rangle:D\in{\cal D}$.

Let $\cal A$ be the family of reflection hyperplanes in $D^n$
for all reflections from $\Gamma$. In other words, $\cal A$ is the
family of all images through elements of $\Gamma$ of the hypeplanes
$A_D:D\in{\cal D}$. Each $A\in{\cal A}$ splits $D^n$ into two components.
Denote the closures of these components in $D^n$ by $Y$ and $Z$,
and observe that $(A,\{ Y,Z \})$ is an elementary splitting of $D^n$.
Denote by $\cal C$ the set of elementary splittings  $(A,\{ Y,Z \})$
as above, for all $A\in{\cal A}$. It is fairy clear that $\cal C$
is then a discrete family of pairwise noncrossing splittings of $D^n$,
i.e. a tree decomposition of $D^n$, and that it is fine.

It remains to show that the tree system $\Theta_{\cal C}$ associated 
to $\cal C$ is isomorphic to ${\cal M}_\partial(D^n)$.
To see this, note that for each $D\in{\cal D}$ the domain $\Omega_{A_D,H_D^c}$
coincides with the fundamental domain $\Omega_0$ for $\Gamma$.
It is not hard too see that this domain is homeomorphic to
the $n$-disk. The separators of $\cal C$ contained in $\Omega_0$
are exactly $A_D:D\in{\cal D}$, and they clearly form a null and
dense family of collared and pirwise disjoint $(n-1)$-disks 
in the boundary $\partial\Omega_0$. Consequently, the constituent space
$\Omega_0$ of $\Theta_{\cal C}$, together with the family
$A_D:D\in{\cal D}$ of its all peripheral subspaces, is homeomorphic
to $(D^n,{\cal D})$. To see that the same is true for all other
constituent spaces of $\Theta_{\cal C}$, note that each such space has a form
$\gamma\Omega_0$ for some $\gamma\in\Gamma$, and the corresponding 
family
of peripheral subspaces has a form $\gamma A_D:D\in{\cal D}$.
This shows that $\Theta_{\cal C}$ is isomorphic to ${\cal M}_\partial(D^n)$,
thus completing the proof.

\medskip\noindent
{\bf 3.B.12.2 Remark.}
Let $\cal C$ be a family of elementary splittings of $D^n$, as constructed
in the proof of Lemma 3.B.12.1. Restricting $\cal C$ to the boundary 
$\partial D^n$, we obviously get a fine tree decomposition of the
$(n-1)$-sphere, and we denote it by ${\cal C}|_{\partial D^n}$.
Moreover, the domains of this new decomposition are 
the intersections of the domains of $\cal C$ with $\partial D^n$,
and it is not hard to verify that they are the $(n-2)$-dimensional
Sierpi\'nski compacta $(S^{n-1})^\circ$. Thus, these Sierpi\'nski compacta
are the constituent spaces of the associated tree system 
$\Theta_{{\cal C}|_{\partial D^n}}$, and the peripheral subspaces
correspond exactly to the peripheral spheres in these Sierpi\'nski
compacta. It follows that the associated tree system 
$\Theta_{{\cal C}|_{\partial D^n}}$ is isomorphic to the tree system
${\cal M}(S^{n-1})$. Since, by Theorem 3.B.10, we have
$$
\lim\Theta_{{\cal C}|_{\partial D^n}}=\partial D^n=S^{n-1},
$$
it follows that $\lim{\cal M}(S^{n-1})=S^{n-1}$.
Thus, we get an alternative (and more elementary) proof of 
Corollary 1.E.3.2.

\medskip\noindent
{\bf 3.B.13 Remark/Example/Exercise.}
Using a similar argumant as in the proof of Lemma 3.B.12.1 one can
identify, up to homeomorphism, limits of various other tree systems.
For example, one can show that the limit of any dense tree system
of internally punctured $n$-disks (see Example 2.D.3) is homeomorphic
to the $(n-1)$-dimensional Sierpi\'nski compactum $(S^n)^\circ$.
Once this is known, one can use a consolidation procedure from Section
3.A to show that the limit of any dense tree system of internally
punctured connected planar surfaces is homeomorphic to the Sierpi\'nski
curve.

\bigskip\noindent
{\bf 3.C Subdivision of a tree system.}

\medskip
Generalizing the concepts from the previous section, we describe
in this section the operation of subdivision of a tree system,
opposite to the operation of consolidation described in Section 3.A.

Let $\Theta=(T,\{K_t\},\{\Sigma_e\},\{\phi_e\})$ be a tree system,
let $t\in V_T$ be any vertex,
and suppose that $\cal C$ is a tree decomposition of the space $K_t$.
We say that $\cal C$ {\it does not cross $\Theta$}
if for each separator $A$ of $\cal C$ and each $e\in N_t$ 
there is a halfspace $H$ for $A$ such that
$\Sigma_e\subset\dot{H}$. This means in particular that
$A\cap\Sigma_e=\emptyset$.

Furthermore, given a tree decomposition $\cal C$
of $K_t$ not crossing $\Theta$, and any edge $e\in N_t$,
we say that $\cal C$ is {\it discrete at $e$} if for some (and hence any)
separator $A$ of $\cal C$ there are only finitely many separators in $\cal C$
that separate $A$ from $\Sigma_e$ (i.e. finitely many separators 
$A'$ from $\cal C$ such
that for some halfspace $H$ related to $A'$ we have $A\subset\dot{H}$
and $\Sigma_e\subset\dot{H}^c$). 
It is not hard to realize that $\cal C$ is discrete at $e$ iff
$\Sigma_e\subset\Omega$ for some domain $\Omega\subset K_t$
for $\cal C$.
We also have the following sufficient condition for discreteness at $e$.

\medskip\noindent
{\bf 3.C.1 Lemma.}
{\it Let $\Theta=(T,\{K_t\},\{\Sigma_e\},\{\phi_e\})$ 
be a tree system, let $\cal C$ be a tree decomposition
not crossing $\Theta$ of some constituent space $K_t$, and let
$e\in N_t$. If $\cal C$ is fine and if $\hbox{diam}(\Sigma_e)>0$,
then $\cal C$ is discrete at $e$.}

\medskip\noindent
{\bf Proof:} Let $A$ be any separator from $\cal C$, and let $H_A$
be this halfspace for $A$ which does not contain $\Sigma_e$.
Let $A'$ be any separator from $\cal C$ that separates $A$ from
$\Sigma_e$, with $A\subset \dot{H}$ and $\Sigma_e\subset\dot{H}^c$
for the corresponding halfspaces $H,H^c$ for $A'$.
Note that $H_A\subset\dot{H}$, and consequently
$$\hbox{diam}(H)\ge\hbox{diam}(H_A)>0
\hbox{\quad and\quad} 
\hbox{diam}(H^c)\ge\hbox{diam}(\Sigma_e)>0.
$$
In view of fineness of $\cal C$, this implies that there are only
finitely many separators $A'$ as above, which completes the proof.

\medskip\noindent
{\bf 3.C.2 Definition.}
A tree decomposition $\cal C$ of a space $K_t$ from
a tree system $\Theta$ is {\it compatible with $\Theta$} if it does
not cross $\Theta$ and if it is discrete at $e$ for each $e\in N_t$.

\medskip\noindent
{\bf 3.C.3 Example.}
Given a tree system $\Theta=(T,\{K_t\},\{\Sigma_e\},\{\phi_e\})$
and a vertex $t\in V_T$, we 
present a criterion for a family of elementary splittings of the space $K_t$ 
to be a tree decomposition of $K_t$ compatible
with $\Theta$. We will use this criterion in Section 3.D.

Let $x_0\in K_t\setminus\cup\{ \Sigma_e:e\in N_t \}$.
Suppose that a family ${\cal C}_t$ of elementary splittings of $K_t$ has a form
$$
{\cal C}_t=\{ (A_k, \{ H_k,H_k^c \}):k\in I\hskip-3ptN \},
$$
where the following three conditions hold:
\itemitem{(1)} splittings from ${\cal C}_t$ do not cross $\Theta$, i.e. each separator $A_k$ is disjoint with each of the sets from the family ${\cal A}_t=\{ \Sigma_e:e\in N_t \}$ and each halfspace $H_k$ is ${\cal A}_t$-saturated,
\itemitem{(2)} $H_{k+1}\subset\dot{H}_k$ for each natural $k$, 
\itemitem{(3)} $\cap_{k=1}^\infty H_k=\{ x_0 \}$ (equivelently, $x_0\in\cap_{k=1}^\infty H_k$ and 
$\lim_{k\to\infty}\hbox{diam}(H_k)=0$).

\noindent
We then get the following.

\medskip\noindent
{\bf 3.C.3.1 Lemma.} {\it If a family ${\cal C}_t$ satisfies
conditions (1)-(3) above then
${\cal C}_t$ is fine,
discrete and discrete at all $e\in N_t$. In particular, it is a tree 
decomposition of $K_t$, and it is compatible with $\Theta$.}

\medskip\noindent
{\bf Proof:} By condition (2), splittings in ${\cal C}_t$ are pairwise noncrossing, and ${\cal C}_t$ is discrete. By condition (3), ${\cal C}_t$ is fine. Thus, in view of condition (1), ${\cal C}_t$ is a tree decomposition of $K_t$. Moreover, by condition (3), for each $e\in N_t$ there is $k$ such that $\Sigma_e\cap H_k=\emptyset$. It follows that for each $e\in N_t$ the family ${\cal C}_t$ is discrete at $e$. This completes the proof.

\medskip\noindent
{\bf 3.C.4 Definition.}
A {\it tree decomposition of a tree system $\Theta$}
is a family ${\cal C}=\{ {\cal C}_t:t\in V_T \}$ 
of tree decompositions ${\cal C}_t$ of the spaces $K_t$
which are all compatible with $\Theta$.

\medskip
We now describe some elementary splittings of the limit space
$\lim\Theta$ induced by any elementary splittings from 
any tree decomposition ${\cal C}=\{ {\cal C}_t \}$ of
$\Theta$. Let $t$ be an arbitrary vertex in the underlying tree $T$
of $\Theta$ and let $(A,\{H,H^c\})$ be an elementary
splitting of the space $K_t$ belonging to ${\cal C}_t$.
Note that both subsets $H$ and $H^c$ of $K_t$ are
saturated with respect to the family ${\cal A}_t=\{ \Sigma_e:e\in N_t \}$.
In particular, it makes sense to speak of the subsets $G(H)$ and $G(H^c)$
in $\lim\Theta$ of the form described in Section 1.C 
(just before Proposition 1.C.1). Moreover, when we view $A$ as the subset
of $\lim\Theta$, we clearly have $G(H)\cap G(H^c)=A$, where all three
subsets in this expression are compact subspaces of $\lim\Theta$.
Thus, the triple $(A,\{ G(H),G(H^c) \})$ is an elementary splitting
of the limit space $\lim\Theta$.
For each $t\in V_T$ we put
$$
{\cal C}_t^{\lim}=\{ (A,\{ G(H),G(H^c) \}): 
A \hbox{ is a separator in } {\cal C}_t \}
$$

\medskip\noindent
{\bf 3.C.5 Lemma.}
{\it Let $\Theta$ be a tree system of metric compacta, and let
${\cal C}(\Theta)$ be the associated tree decomposition of the limit $\lim\Theta$.
Then for any tree decomposition ${\cal C}$ of $\Theta$ the family
$$
{\cal C}_{\lim}=(\bigcup_{t\in V_T}{\cal C}_t^{\lim})
\cup{\cal C}(\Theta)
$$ 
is a tree decomposition of $\lim\Theta$.}

\medskip\noindent
{\bf Proof:} The fact that the elementary splittings from ${\cal C}_{\lim}$
are pairwise noncrossing follows directly from the fact 
that the families ${\cal C}_t$
do not cross $\Theta$. The discreteness of ${\cal C}_{\lim}$ follows
fairly directly from discreteness of each ${\cal C}_t$ at each $e\in N_t$.
We omit further details.

\medskip\noindent
{\bf 3.C.6 Proposition.} {\it Under notation of Lemma 3.C.5,
the tree decomposition ${\cal C}_{\lim}$ is fine iff each tree decomposition 
${\cal C}_t$ from $\cal C$ is fine.}

\medskip\noindent
{\bf Proof:}
One implication,
namely that fineness of ${\cal C}_{\lim}$ implies fineness of every 
${\cal C}_t$, is immediate just by observing that halfspaces for ${\cal C}_t$
are simply restrictions to $K_t$ of the appropriate halfspaces for
${\cal C}_{\lim}$. The converse implication 
requires much more effort and some preparatory claims.
We start with a claim which gives a useful characterization
of fineness of a tree decomposition. 

\medskip\noindent
{\bf Claim 1.} {\it A tree decomposition $\cal C$ of a compact $K$
is fine iff for any $\epsilon>0$ there is a finite collection
$\Omega_1,\dots,\Omega_m$ of domains for $\cal C$
such that if $H$ is any halfspace from $\cal C$ not containing any
of the above domains then $\hbox{diam}(H)<\epsilon$.}

\medskip
To prove Claim 1, suppose first that $\cal C$ is fine. Given $\epsilon>0$,
let $(A_i,\{ Y_i,Z_i \}):i=1,\dots,q$ be all elementary splittings in $\cal C$
for which $\hbox{diam}(Y_i)\ge\epsilon$ and 
$\hbox{diam}(Z_i)\ge\epsilon$. Let $\Omega_1,\dots,\Omega_m$
be the set of all domains for $\cal C$ adjacent to some of the separators
$A_1,\dots,A_q$. Since each separator has exactly two 
adjacent domains, the above set of domains is finite.
Let $H$ be a halfspace from $\cal C$ not containing any of the
above domains. It is not hard to realize that then for any $1\le i\le q$
either $Y_i$ or $Z_i$ is contained in $H^c$, and hence
$\hbox{diam}(H^c)\ge\epsilon$. On the other hand, the separator $A$
corresponding to $H$ is clearly distinct from each of the separators
$A_1,\dots,A_q$, and hence $\hbox{diam}(H)<\epsilon$,
as required.

To prove the converse implication in Claim 1, fix $\epsilon>0$
and let $\Omega_1,\dots,\Omega_m$ be some domains associated
to $\epsilon$, as in the assumtion of the implication.
Note that, by the assumtion, for any splitting $(A,\{ Y,Z \})\in{\cal C}$ 
with at least one halfspace not containing any of the above domains
we have
$\min(\hbox{diam}(Y),\hbox{diam}(Z))<\epsilon$.
Thus, it is sufficient to show that the number of splittings with
both halfspaces containing some of the domains 
$\Omega_1,\dots,\Omega_m$ is finite. Translating this to the
language of the dual tree, we need to know that for any finite
set $V_0$ of vertices in $T_{\cal C}$ the set of edges in $T_{\cal C}$
which separate some two of the vertices of $V_0$ is finite.
Since this set of edges clearly coincides with the set of edges 
in the subtree of $T_{\cal C}$ spanned by $V_0$, this completes
the proof of Claim 1.  

\medskip
We come back to the tree decomposition ${\cal C}_{\lim}$.
For each $t\in V_T$ and any $e\in N_t$ let $H_e$ be this halfspace
from ${\cal C}(\Theta)$ associated to the separator $\Sigma_e$
which does not contain $K_t$. For each halfspace $H$ from ${\cal C}_t$
put ${\cal H}_H:=\{ H_e:e\in N_t \hbox{ and } \Sigma_e\subset H \}$.
The next preparatory claim provides some estimate for the diameter
of a halfspace from ${\cal C}_{\lim}$ in terms of diameters
of appropriate halfspaces from $\cal C$ and ${\cal C}(\Theta)$. 

\medskip\noindent
{\bf Claim 2.}
{\it Let $t$ be any vertex of $T$ and let $H$ be any halfspace from
${\cal C}_t$. Then the induced halfspace $G(H)$ from ${\cal C}_t^{\lim}$
has form $G(H)=H\cup(\bigcup{\cal H}_H)$ and its diameter is estimated by
$$
\hbox{diam}(G(H))\le\hbox{diam}(H)+2\cdot
\max\{ \hbox{diam}(H'):H'\in{\cal H}_H \}.
$$}

The first assertion of Claim 2, i.e. that $G(H)=H\cup(\bigcup{\cal H}_H)$,
follows directly from the definition of $G(H)$.
To prove the second assertion, we estimate the distance of any two points
of $G(H)$. Suppose that $x,y$ are some points of $G(H)$
not contained in $H$. Then, by the first assertion, there exist $e,e'\in N_t$
such that $x\in H_e$ and $y\in H_{e'}$.
Let $x'\in\Sigma_e$ and $y'\in\Sigma_{e'}$ be arbitrary points.  
Since $\Sigma_e=H\cap H_e$ and $\Sigma_{e'}=H\cap H_{e'}$,
we get the estimate
$$
d(x,y)\le d(x,x')+d(x',y')+d(y',y)\le\hbox{diam}(H_e)+\hbox{diam}(H)
+\hbox{diam}(H_{e'})
$$
which can be further estimated from above by the right hand side
term in the desired inequality. In case when one or both point $x,y$
belong to $H$ a similar (even simpler) estimate can be obtained
in the same way. This completes the proof of Claim 2.

\medskip
We come back to the proof of the following implication:
if each tree decomposition ${\cal C}_t$ from $\cal C$ is fine
then ${\cal C}_{\lim}$ is fine. We use the characterization of
fineness given in Claim 1. Fix any $\epsilon>0$.
Since ${\cal C}(\Theta)$ is fine (see the last comment in Example
3.B.7),
we choose a finite subtree $S\subset T$ such that any
halfspace $H'$ from ${\cal C}(\Theta)$ which does not contain
any of the spaces $K_t:t\in V_S$ satisfies 
$\hbox{diam}(H')<{\epsilon\over3}$.
Note that this implies also that for any $s\in V_T\setminus V_S$
we have $\hbox{diam}(K_s)<{\epsilon\over3}$.
For each $t\in V_S$, using the fact that ${\cal C}_t$ is fine,
choose some domains $\Omega^t_1,\dots,\Omega^t_{m_t}$,
with $m_t\ge1$, such that
any halfspace $H$ from ${\cal C}_t$ which does not contain
any of these domains satisfies $\hbox{diam}(H)<{\epsilon\over3}$.
Put ${\cal Z}:=\cup_{t\in V_S}\{ \Omega^t_1,\dots,\Omega^t_{m_t} \}$
and note that $\cal Z$ is a finite family of domains for the tree
decomposition ${\cal C}_{\lim}$.
We claim that any halfspace $H_{\lim}$ from ${\cal C}_{\lim}$ which does not
contain any of the domains from $\cal Z$ satisfies
$\hbox{diam}(H_{\lim})<\epsilon$.

To prove the above claim, suppose first that $H_{\lim}=H$
for some halfspace $H$ from ${\cal C}(\Theta)$. Since $H=H_{\lim}$
does not contain any of the domains from $\cal Z$, it also does not
contain any of the spaces $K_t:t\in V_S$. By the choice of $S$,
this implies that $\hbox{diam}(H_{\lim})<{\epsilon\over3}<\epsilon$,
as required. 

In the remaining case we have $H_{\lim}=G(H)$, where $H$
is a halfspace from ${\cal C}_t$ in the space $K_t$, for
some $t\in V_T$. We will estimate the diameter of $H_{\lim}=G(H)$
using Claim 2. To do this, we claim that 
$\hbox{diam}(H)<{\epsilon\over3}$. Indeed, if $H\subset K_t$
for some $t\in V_S$, this estimate follows from our choice
of the domains $\Omega^t_1,\dots,\Omega^t_{m_t}$ in view
of the fact that $H$ does not contain any of them.
If $H\subset K_s$ for some $s\in V_T\setminus V_S$,
we get $\hbox{diam}(H)\le\hbox{diam}(K_s)<{\epsilon\over3}$,
where the last inequality follows from the choice of $S$.

We now estimate diameters of the halfspaces $H'\in{\cal H}_H$.
Since any such $H'$ (being a halfspace from ${\cal C}(\Theta)$)
does not contain any of the spaces $K_t:t\in V_S$
(because it does not contain any of the domains from $\cal Z$),
it follows from the choice of $S$ that 
$\hbox{diam}(H')<{\epsilon\over3}$.

In view of Claim 2, the estimates from the two previous paragraphs
yield the inequality $\hbox{diam}(G(H))<\epsilon$, as required.
By Claim 1, the tree decomposition ${\cal C}_{\lim}$ is then fine,
which completes the proof.

\medskip\noindent
{\bf 3.C.7 Definition.}
Let $\Theta$ be a tree system of metric compacta.
A {\it subdivision} of $\Theta$ 
is any tree system of form $\Theta_{{\cal C}_{\lim}}$, for any fine tree
decomposition ${\cal C}_{\lim}$ as in Lemma 3.C.5.

\medskip\noindent
{\bf 3.C.8 Proposition.} {\it Let $\Xi$ be any subdivision of $\Theta$.
Then $\Theta$ can be canonically obtained from $\Xi$ by means of 
a consolidation. Moreover, the limits $\lim\Xi$ and $\lim\Theta$
are canonically homeomorphic.}

\medskip\noindent
{\bf Proof:}
Let $\Xi=\Theta_{{\cal C}_{\lim}}$.
Note that the constituent spaces of the tree system $\Xi$
are precisely the constituent spaces of the systems $\Theta_{{\cal C}_t}$
for all $t\in V_T$.
Since each ${\cal C}_t$ is a fine tree decomposition of $K_t$, it follows
from Theorem 3.B.10 that $\lim\Theta_{{\cal C}_t}=K_t$.
This shows that $\Theta$ is a consolidation of $\Xi$, for the canonical
partition of the dual tree $T_{{\cal C}_{\lim}}$ into subtrees $T_{{\cal C}_t}$.
The second assertion follows either from Theorem 3.B.10 and Lemma 3.C.5
or, in view of the first assertion, from Theorem 3.A.1.


\bigskip\noindent
{\bf 3.D Homogeneity of trees of manifolds.}

\medskip
In this section we prove the following.

\medskip\noindent
{\bf 3.D.1 Proposition.}
{\it Let $M$ be a closed connected topological manifold
(either oriented or non-orientable).
Then the tree of manifolds $M$, i.e. the Jakobsche space ${\cal X}(M)$
(as defined in Section 1.E), is homogeneous.}

\medskip
This result has been proved before in [J2] (for oriented $M$)
and in [St] (for non-orientable $M$ which are PL). 
The proof we present here uses
the technique of subdivisions and consolidations of tree systems.
It is inspired by the corresponding proof in Sections 7 and 8 of [J2].
Our argument,
and the technique used in it, has a potential for extensions.
It can be applied
to various other classes of tree systems (e.g. to tree systems of polyhedra
mentioned in the introduction), to study orbits of the group
of homeomorphisms of the corresponding limit space.
One easy instance of such extension is presented below, as
Propositions 3.D.6.1 and 3.D.6.2. 

We start with few technical preparatory results.

\medskip\noindent
{\bf 3.D.2 Lemma.}
{\it Let $M$ be a closed connected topological manifold, oriented or non-orientable,
and let $\cal M$ be the dense tree system of manifolds $M$, 
with the underlying tree $T$. Then the points
of the limit ${\cal X}(M)=\lim{\cal M}$ corresponding to the set $E_T$
of the ends of $T$ are all in the same orbit of the group of
homeomorphisms of ${\cal X}(M)$.}

\medskip\noindent
{\bf Proof:} Let $z_1,z_2\in E_T\subset\lim{\cal M}$. 
Using Toru\'nczyk's Lemma 1.E.2.1
(together with Lemma 1.E.4.1 in case when $M$ is non-orientable)
it is not hard to get an automorphism of the tree system $\cal M$
for which the corresponding automorphism of the underlying tree $T$
maps $z_1$ to $z_2$. Clearly, this automorphism of $\cal M$
induces a homeomorphism of $\lim{\cal M}$ which maps $z_1$ to $z_2$.
This finishes the proof.

\medskip
Next result is an extension of Lemma 5 from [J1], and it appears
implicitely inside the proof of Lemma 7.1 in [J2].

\medskip
We will need the following corollary to Lemma 1.D.2.2.

\medskip\noindent
{\bf 3.D.3 Corollary.}
{\it Let $M$ and $\cal D$ be as in Lemma 1.D.2.2, and let $x_0$
be an interior point of $M$ not contained in any $D\in{\cal D}$.
Then there is a sequence of collared $n$-disks $Q_k$ contained 
in $\hbox{int}(M)$, such that}
\item{(1)} {\it for each $k$ the boundary $\partial Q_k$ is disjoint with 
the union of $\cal D$,}
\item{(2)} {\it for each $k$ we have $Q_{k+1}\subset\hbox{int}(Q_k)$,}
\item{(3)} {\it $\{ x_0 \}=\cap_k Q_k$.}

\medskip\noindent
{\bf Proof:}
By Lemma 1.D.2.2,
the quotient $M/{\cal D}$ is homeomorphic to $M$.
Let $A_{\cal D}$ be the set of points in the quotient $M/{\cal D}$
corresponding to the collapsed elements of $\cal D$.  Clearly, 
$A_{\cal D}$ is then countable infinite, and 
$x_0\in (M/{\cal D})\setminus A_{\cal D}$. Obviously, there exists
a sequence of collared $n$-disks $P_k$ in $M/{\cal D}$ such that
\item{(1)} for each $k$ the boundary $\partial P_k$
is disjoint with $A_{\cal D}$,
\item{(2)} for each $k$ we have $P_{k+1}\subset\hbox{int}(P_k)$,
\item{(3)} $x_0=\cap_k P_k$.

\noindent
For each $k$ put $Q_k=q^{-1}(P_k)$, where $q:M\to M/{\cal D}$ 
is the quotient map. To see that the sequence $Q_k$ is as required,
it obviously suffices to show that each $Q_k$ is a collared $n$-disk
in $M$. To do this, consider the decomposition of $Q_k$ induced
by the family ${\cal D}_k=\{ D\in{\cal D}:D\subset Q_k \}$.
Then we clearly have $P_k=Q_k/{\cal D}_k$, and we denote the
quotient map by $q_k$. Arguing as in the proof of Lemma 1.D.2.2,
we get that the above decomposition is shrinkable, and hence
$q_k$ can be approximated by homeomorphisms.
This shows $Q_k$ is an $n$-disk. By applying the same argument
to the complement $M\setminus\hbox{int}(Q_k)$ (instead of $Q_k$),
we get that $Q_k$ is collared, which finishes the proof.

\medskip
A crucial ingredient in the proof of Proposition 3.D.1,
is the following result, in the proof of which we use the 
technique of modifications of tree systems, developed 
in Sections 3.A-3.C above. This proof is inspired by the Jakobsche's idea
from his proof of Lemma (7.1) in [J2].

\medskip\noindent
{\bf 3.D.4 Lemma.}
{\it Let $M$ be a closed connected manifold, oriented or non-orientable,
and let ${\cal M}=(T,\{ K_t \},\{ \Sigma_e \},\{ \phi_e \})$ 
be the dense tree system of manifolds $M$. Let 
$$
x\in K_{t_*}\setminus\cup\{ \Sigma_e:e\in N_{t_*} \}\subset\lim{\cal M}
$$
for some $t_*\in V_T$,
and let $y\in E_T\subset\lim{\cal M}$.
Then there is a homeomorphism of $\lim{\cal M}$ which maps $x$ to $y$.}

\medskip\noindent
{\bf Proof:}
Denoting $n=\dim M$, 
recall that for all $t\in V_T$ we have
$$
K_t=M_t^\circ=M_t\setminus\cup\{ \hbox{int}(D):D\in{\cal D} \},
$$
where $M_t$ is a homeomorphic copy of $M$, and
where $\cal D$ is a null and dense family of pairwise disjoint
collared $n$-disks in $M_t$. View $x$ as a point of $M_{t_*}$ and
consider a sequence $Q_k$ of $n$-disks as in Corollary 3.D.3,
for $x_0=x$. For each $k$ consider the elementary splitting 
$(A_k,\{ H_k,H_k^c \})$ of $K_{t_*}$ given by
$A_k:=\partial Q_k$ and $H_k:=Q_k\cap K_{t_*}$, 
and denote by ${\cal C}_{t_*}$
the set of these elementary splittings for all natural $k$. 
By Lemma 3.C.3.1, the family ${\cal C}_{t_*}$ is 
a tree decomposition of $K_{t_*}$ compatible with $\cal M$.

Let ${\cal C}^x=\{ {\cal C}_t:t\in V_T \}$ be a tree
decomposition of $\cal M$ such that 
${\cal C}_{t_*}$ is as above and
${\cal C}_t=\emptyset$ for $t\ne t_*$.
Let ${\cal M}^x=(T^x,\{ K_t^x \},\{ \Sigma_e^x \},\{ \phi_e^x \})$ 
be the subdivision of $\cal M$ induced by ${\cal C}^x$. 
The underlying tree $T^x$ may be viewed as obtained from $T$
by expanding the vertex $t_*$ into an infinite polygonal ray 
$\rho=(t_0,t_1,\dots)$. We may identify all other vertices of $T^x$
bijectively with the vertices in $V_T\setminus\{ t_* \}$.
The set of the edges of $T$ adjacent to $t_*$ canonically splits
into subsets which may be identified with the sets of edges of $T^x$
adjacent to the vertices $t_i$ (for $i=0,1,\dots$). 
The edges of $T^x$ disjoint
from the ray $\rho$ are in the natural bijective corresponedence
with the deges of $T$ not adjacent to $t_*$. 
Accordingly, the constituent spaces $K_t:t\ne t_*$ are still the
constituent spaces in ${\cal M}^x$ at the corresponding vertices,
with the families of peripheral subspaces unchanged.
In particular, they are still the densely punctured manifolds $M$.
The constituent space $K_{t_*}$ of $\cal M$ splits into the 
constituent spaces $K^x_{t_i}:i\ge0$ of ${\cal M}^x$,
which are described as follows.
We have 
$K^x_{t_0}:=H_1^c\cap K_{t_*}=K_{t_*}\setminus\hbox{int}(Q_1)$,
and for $i\ge1$ we have 
$K^x_{t_i}:=(H_{i+1}^c\setminus\dot{H}_i)\cap K_{t_*}=K_{t_*}\cap(Q_i\setminus\hbox{int}(Q_{i+1}))$.
It is not hard to note that
$K^x_{t_0}$ still has the form of the denesely punctured manifold $M$,
while each of $K^x_{t_i}:i\ge1$
is the densely punctured sphere $S^n$.

By Proposition 3.C.8, we have the canonical
identification of the limits $\lim{\cal M}$ and $\lim{\cal M}^x$.
The point $x$, viewed as an element of $\lim{\cal M}^x$, clearly
corresponds to the end  
of the tree $T^x$ induced by the ray $\rho$,
i.e. $x=[\rho]\in E_{T^x}$.
On the other hand, the point $y$ viewed as an element of $\lim{\cal M}^x$ 
still corresponds to an end
of the underlying tree, i.e $y\in E_{T^x}\subset\lim{\cal M}^x$.

We now apply to the system ${\cal M}^x$ an operation of consolidation,
as described in Section 3.A. More precisely, for each $i\ge1$ choose any 
vertex $s_i$ in $T^x$ adjacent to $t_i$ and distinct from both $t_{i-1}$
and $t_{i+1}$, and denote by $\Pi$ the partition of $T^x$ consisting
of the subtrees $S_i:i\ge1$ spanned on the pairs $t_i,s_i$
(these subtrees are just edges), and of subtrees reduced to vertices
for all remaining vertices of $T^x$. 
Let ${\cal M}^x_\Pi$ be the tree system obtained from ${\cal M}^x$
by consolidation with respect to $\Pi$. By Proposition 3.C.8,
the limit $\lim{\cal M}^x_\Pi$ is canonically homeomorphic to
$\lim{\cal M}^x$, and hence also to $\lim{\cal M}$.

From what was said above about the system ${\cal M}^x$
it is not hard to deduce that ${\cal M}^x_\Pi$ is a dense tree system
of manifolds $M$. Moreover, the point $x$ (viewed now as a point
of $\lim{\cal M}^x_\Pi$) clearly corresponds to the end of the tree
$T^x_\Pi$ induced by the ray $\rho_\Pi=(S_1,S_2,\dots)$.
On the other hand, the point $y$
still corresponds to an end of the underlying tree, i.e. $y\in E_{T^x_\Pi}\subset \lim{\cal M}^x_\Pi$. Thus, by Lemma 3.D.2,
there is a homeomorphism of $\lim{\cal M}^x_\Pi$ which
maps $x$ to $y$, which finishes the proof.

\medskip
The next result is a consequence of Lemma 3.D.4.

\medskip\noindent
{\bf 3.D.5 Corollary.}
{\it Let $M$ and $\cal M$ be as in Lemma 3.D.4,
let $x\in\Sigma_{e_*}\subset\lim{\cal M}$ for some $e_*\in O_T$,
and let $y\in E_T\subset\lim{\cal M}$.
Then there is a homeomorphism of $\lim{\cal M}$ which maps $x$ to $y$.}

\medskip\noindent
{\bf Proof:} Denoting $n=\dim M$,
view the constituent space $K_{\alpha(e_*)}$ of $\cal M$ as obtained
from $M_{\alpha(e_*)}=M$ by deleting interiors of disks $D$ from
a null and dense family $\cal D$ of pairwise disjoint collared $n$-disks.
Let $D_*\in{\cal D}$ be this disk for which $\partial D_*=\Sigma_{e_*}$.
Let $Q$ be a collared $n$-disk in $M_{\alpha(e_*)}$ with
$D_*\subset Q$ and with $\partial Q$ disjoint from the union of $\cal D$
(existence of such $Q$ follows by the arguments as in the proof of
Lemma 3.D.3).

Consider the elementary splitting $(A,\{ H,H^c \})$ of $K_{\alpha(e_*)}$
given by $A=\partial Q$ and $H=Q\cap K_{\alpha(e_*)}$.
Let ${\cal M}^Q$ be the subdivision of $\cal M$ induced by this
single splitting. Let $t_H$ be the vertex of the underlying tree $T^Q$
of ${\cal M}^Q$ corresponding to the constituent space $H$.
Viewing $\omega(e_*)$ naturally as a vertex in $T^Q$, note that 
$t_H$ and $\omega(e_*)$ are adjacent, and denote by $S_Q$
the subtree of $T^Q$ consisting of these two vertices and the edge
which connects them. Consider the consolidation 
${\cal M}^Q_\Pi$ of ${\cal M}^Q$ for the partition $\Pi$ consisting
of the subtree $S_Q$ and the singleton subtrees for all other vertices.
Clearly, the limits $\lim{\cal M}^Q_\Pi$ and $\lim{\cal M}$ are
canonically homeomorphic.

Note that ${\cal M}^Q_\Pi$ is again a dense tree of manifolds $M$,
and that $x$, naturally viewed as element of $\lim{\cal M}^Q_\Pi$,
belongs to the constituent space $K_{S_Q}\subset\lim{\cal M}^Q_\Pi$,
and lies outside its all peripheral subspaces.
Since the point $y$, viewed as element of $\lim{\cal M}^Q_\Pi$,
still corresponds to an end of the underlying tree, 
the corollary follows by applying Lemma 3.D.4.

\medskip\noindent
{\bf Proof of Proposition 3.D.1:}
The proposition is a direct consequence of Lemmas 3.D.2, 3.D.4
and Corollary 3.D.5.


\medskip\noindent
{\bf 3.D.6 Remarks.}

The argument as above in this section yields also the following.

\medskip\noindent
{\bf 3.D.6.1 Proposition.}
{\it Let $M$ be a connected compact topological manifold with boundary,
either oriented or non-orientable, and let $\cal M$
be the dense tree system of internally punctured manifolds $M$,
as defined in Remark 2.D.3.1.
Then, viewing the boundaries $\partial M_t$ 
as subsets of the
constituent spaces $K_t=M_t^\circ$ of $\cal M$,
all points of $K_t\setminus\partial M_t\subset\lim{\cal M}$ 
(for all $t$) and all points of $E_T\subset\lim{\cal M}$
are in the same orbit of the group of homeomorphisms of 
$\lim{\cal M}={\cal X}_{int}(M)$.}

\medskip
Recall that, given
any natural number $m$,  
a topological space $X$ is {\it $m$-homogeneous} if the group of its
homeomorphisms acts transitively on the set of all $m$-tuples
of pairwise distinct points of $X$.
A straightforward extension of the arguments
of this section allows to prove Theorem 8.1 of [J2], i.e.
the fact that the Jakobsche space
${\cal X}(M)$ is $m$-homogeneous, for arbitrary $m$.
It also allows to prove the following variant of this result,
for trees of manifolds with boundary.

\medskip\noindent
{\bf 3.D.6.2 Proposition.} 
{\it Let $M$ and $\cal M$
be as in Proposition 3.D.6.1.
Consider the subspace 
$$
U=E_T\cup(\bigcup_{t\in V_T}K_t\setminus\partial M_t)\subset\lim{\cal M}.
$$
Then for any natural number $m$ the group of homeomorphisms of $\lim{\cal M}$
acts transitively on the set of $m$-tuples of pairwise distinct points of $U$.}


\bigskip\noindent
{\bf 3.E Weakly saturated tree systems of manifolds.}

\medskip
In this section we present another application of the operations of consolidation and subdivision. It concerns dense trees of finite
families of manifolds, as in Example 3.A.2, and provides a significant
strengthening of Propositions 3.A.2.1 and 3.A.2.3. More precisely, we will show
that under a much weaker assumption than 2-saturation, the limit of a dense
tree of a finite family of manifolds $\{ M_1,\dots,M_k \}$ is still homeomorphic to the space ${\cal X}(M_1\#\dots\# M_k)$. The weaker condition will be called 
{\it weak saturation}. The results of this section are used by the author 
(in another paper [Sw2]) to show that various trees of manifolds,
in arbitrary dimension, appear as Gromov boundaries of some hyperbolic
groups.

\medskip
Let ${\cal N}=\{ M_1,\dots,M_k \}$ be a finite family of closed
connected topological manifolds of the same dimension, either all oriented,
or at least one of which is non-orientable. Let
${\cal M}=(T,\{ K_t\},\{ \Sigma_e \},\{ \phi_e \})$ be a dense
tree system of manifolds from $\cal N$. As in Example 3.A.2, for each $t\in V_T$
let $i_t\in\{ 1,\dots,k \}$ be this index for which $K_t$
has a form $M_{i_t}^\circ$.
Recall also that {\it a half-tree} in a tree is its any maximal subtree obtained by
deleting the interior of any edge.

\medskip\noindent
{\bf 3.E.1 Definition}
We say that a dense tree system $\cal M$ of manifolds from $\cal N$ is {\it weakly saturated} if for each $j\in\{ 1,\dots,k \}$ any half-tree in the underlying
tree $T$ contains a vertex $t$ with $i_{t}=j$ (equivalently, 
for each $j\in\{ 1,\dots,k \}$ the set
$V^j_T=\{ t\in V_T:i_t=j \}$ spans $T$).

\medskip
The main result of this section is the following.

\medskip\noindent
{\bf 3.E.2 Theorem.}
{\it Let $\cal M$ be any weakly saturated dense tree system of manifolds from a finite family ${\cal N}=\{ M_1,\dots,M_k \}$ 
(where $M_i$ are closed connected topological manifolds of the same
dimension, either all oriented, 
or at least one of which is non-orientable). Then the limit $\lim{\cal M}$ 
is homeomorphic to the Jakobsche space ${\cal X}(M_1\#\dots\# M_k)$.}

\medskip\noindent
{\bf Proof:} The proof consists of showing that, by applying certain 
consolidation followed by some subdivision, the tree system $\cal M$ can be
transformed into a 2-saturated dense tree system of manifolds from $\cal N$.
Due to invariance of the limit under operations as above, and in view of
Propositions 3.A.2.1 and 3.A.2.3, this will give the assertion.

We use the following notation. Given an oriented edge $e\in O_T$ (where $T$
is the underlying tree of the system $\cal M$), denote by $T_e$ this half-tree  obtained by deleting from $T$ the interior of $e$, which contains the 
terminal vertex $\omega(e)$.

We describe appropriate consolidation and subdivision simultaneously,
accordingly with the following scheme. 
We successively choose finite subtrees $S$
of a partition $\Pi$ of the tree $T$ 
(which will induce a desired consolidation ${\cal M}_\Pi$ of $\cal M$),
and finite tree decompositions ${\cal C}_S$ of the corresponding 
constituent spaces $K_S$ in ${\cal M}_\Pi$
(which will give the tree decomposition 
${\cal C}=\{ {\cal C}_S:S\in\Pi \}$  inducing a desired subdivision
of ${\cal M}_\Pi$). The choices of subtrees $S$ and decompositions ${\cal C}_S$
are made inductively, using an auxilliary ordering of the vertices of $T$ into
a sequence $(u_n)_{n\ge1}$, as follows.

\smallskip\noindent
{\it Step 1.} Put $S_1=\{ u_1 \}$ (i.e. the subtree consisting of a single
vertex $u_1$), and put $\Pi_1=\{ S_1 \}$. Furthermore, 
let ${\cal C}_{S_1}$ be the empty tree decomposition
of the corresponding space $K_{S_1}=K_{u_1}$.

\smallskip\noindent
{\it Step 2.}
Having constructed a finite family $\Pi_n$ of finite subtrees, and a corresponding
family of decompositions ${\cal C}_{S}:S\in\Pi_n$, we keep as a part of inductive assumption the following properties (which clearly hold true for $n=1$): 
\item{(i0)} the subtrees in $\Pi_n$ are pairwise disjoint,
\item{(i1)} $u_n\in\cup\{ V_S:S\in\Pi_n \}$,
\item{(i2)} the union $\cup\{ V_S:S\in \Pi_n \}$ is the vertex set of 
a finite subtree of $T$, which we denote $T_n$,
\item{(i3)} for each $S\in\Pi_n$, for each separator $A$ of 
the decomposition ${\cal C}_S$, and for any edge $e\in N_S$ there is 
a halfspace $H$ for $A$ such that $\Sigma_e\subset\dot H$.

\noindent
For each subtree $S\in\Pi_n\setminus\Pi_{n-1}$ and for each domain
$\Omega\subset K_S$ corresponding to ${\cal C}_S$, choose arbitrary 
pairwise distinct oriented edges
$e_1,\dots,e_k,e_1',\dots,e_k'$ from $N_{S}$, not belonging to $T_n$,
and such that for $1\le j\le k$ the corresponding peripheral subsets 
$\Sigma_{e_j}$ and $\Sigma_{e_j'}$ are contained in $\Omega$. 
This is possible since, due to finiteness of ${\cal C}_S$, $\Omega$ 
has nonepty interior in $K_S$, and by denseness of $\cal M$, it thus 
contains infinitely many peripheral subsets $\Sigma_e$ with $e\in N_S$
and not belonging to $T_n$.
For each $1\le j\le k$ choose a vertex $t_j$ in $T_{e_j}$ and $t_j'$ 
in $T_{e_j'}$ such that $i_{t_j}=i_{t_j'}=j$.
This is possible since $\cal M$ is weakly saturated.
Further, for each $j$ choose a finite subtree $S_j(\Omega)$ in $T_{e_j}$ containing
the vertices $\omega(e_j)$ and $t_j$, and  a finite subtree 
$S_j'(\Omega)$ in $T_{e_j'}$ containing the vertices $\omega(e_j')$ and $t_j'$.
We also require that, putting 
$$
\Pi_{n+1}=\Pi_n\cup\bigcup_\Omega\{ S_1(\Omega),\dots,S_k(\Omega),S_1'(\Omega),\dots,S_k'(\Omega) \}
$$
(where $\Omega$ runs through all domains in all spaces $K_S:S\in\Pi_n\setminus\Pi_{n-1}$),
we have $u_{n+1}\in\cup\{ V_S:S\in\Pi_{n+1} \}$.
This clearly holds true if $u_{n+1}\in\cup\{ V_S:S\in\Pi_n \}$; otherwise,
this can be assured as follows. Let $s$ be the vertex in $T_n$ which is 
closest to $u_{n+1}$ in $T$, and let $e$ be the first oriented edge
on the path from $s$ to $u_{n+1}$. Let $S\in\Pi_n$ be the subtree for which
$s\in V_S$, and let $\Omega\subset K_S$ be the domain corresponding to ${\cal C}_S$
for which $\Sigma_e\subset\Omega$. 
We then choose $e_1$ as above so that additionally
we have $u_{n+1}\in T_{e_1}$ (i.e. we put $e_1=e$), 
and then choose $S_1(\Omega)$ so that it contains $u_{n+1}$.
As a consequence of all our choices above, the family $\Pi_{n+1}$
satisfies conditions (i0)-(i2). For example, (i0) holds because if
$S$ and $S'$ are distinct elements of $\Pi_{n+1}\setminus\Pi_n$
then $S\subset T_e$ and $S'\subset T_{e'}$ where $e$ and $e'$
are distinct elements of $N_{T_n}$.

Now, for each $S=S_j(\Omega)$ or $S=S_j'(\Omega)$ as above, we choose
an appropriate tree decomposition ${\cal C}_S$ of the space $K_S$.
To describe it, note that $K_S$ (together with its peripheral subspaces of 
the system ${\cal M}_\Pi$) is homeomorphic to the densely punctured
manifold, denoted $M_S$, which is a
connected sum of the manifolds $M_t:t\in V_S$. 
We denote by $\Delta_e$ the disks in $M_S$ corresponding to the 
peripheral subspaces $\Sigma_e$ of $K_S$. We also denote by $\Delta_j$
and $K_j$ the spaces $\Delta_{\bar e_j}$ and $K_{t_j}$ if $S=S_j(\Omega)$,
and the spaces $\Delta_{\bar e_j'}$ and $K_{t_j'}$ if $S=S_j'(\Omega)$.

Choose any $\Delta_{e_0}\subset M_S$ such that 
$\Sigma_{e_0}\subset K_j$ and note that, by applying Toru\'nczyk's
Lemma 1.E.2.1 to the manifolds $M_S\setminus\hbox{int}(\Delta_{e_0})$
and $M_S\setminus\hbox{int}(\Delta_j)$, we get a homeomorphism 
$h:M_S\to M_S$
(preserving the orientation if all manifolds in $\cal N$ are oriented)
which maps $\Delta_{e_0}$ onto $\Delta_j$, and which preserves the family
of all disks $\Delta_e$ in $M_S$.
We denote by $h^\circ:K_S\to K_S$ the restricted homeomorphism
of the densely punctured manifold.
Now, we consider the finite tree decomposition ${\cal C}_S$ of $K_S$ 
which is induced
by pushing through $h^\circ$ the original tree decomposition of $K_S$
into constituent spaces $K_t:t\in V_S$ (of the system $\cal M$ restricted to $S$).
Obviously, ${\cal C}_S$ satisfies property (i3) above, and it also has
the following property: 

\item{($*$)}
all domains $\Omega\subset K_S$ for ${\cal C}_S$
are densely punctured manifolds from $\cal N$, and the domain $\Omega$
which contains the peripheral subspace $\Sigma_j=\partial\Delta_j$ 
is homeomorphic to $K_j$, i.e. to the densely punctured manifold $M_j$.

\smallskip
We now put $\Pi=\bigcup_{i=1}^\infty\Pi_n$ and note that, by conditions
(i0) and (i1), $\Pi$ describes a decomposition of the tree $T$ into finite subtrees.
We thus consider the induced consolidation ${\cal M}_\Pi$. By codition (i3),
and by finiteness of the decompositions ${\cal C}_S$, the family
${\cal C}=\{ {\cal C}_S:S\in\Pi \}$ is a tree decomposition of the system
${\cal M}_\Pi$. Denoting by ${\cal M}'$ the tree system obtained from
${\cal M}_\Pi$ by the subdivision induced by $\cal C$ (i.e. putting
${\cal M}'=({\cal M}_\Pi)_{{\cal C}_{\rm lim}}$),
we get from the construction, and in particular from the property ($*$) above,
that ${\cal M}'$ is a 2-saturated dense tree system of manifolds from $\cal N$.
By the comment in the first paragraph of the proof, Theorem 3.E.2 follows.

\bigskip
As an application of Theorem 3.E.2, we now describe a class of inverse
sequences of manifolds whose limits are the Jakobsche spaces 
${\cal X}(M_1\#\dots\# M_k)$. This class of inverse sequences is much 
more flexible, and much
more convenient to deal with, than the corresponding class considered
by Jakobsche in [J2] (compare Remark 3.A.2.4). 
For this reason, it can be more efficiently used to identify boundaries
of some spaces and groups as appropriate trees of manifolds, see [Sw2]. 

\medskip\noindent
{\bf 3.E.3 Definiton.}
Let $\cal N$ be a finite family of closed connected $n$-dimensional 
topological manifolds,
either all oriented, or at least one of which is non-orientable.
Let $${\cal J}=(\{ X_i:i\ge1 \},\{ \pi_i:i\ge1 \})$$
be an inverse sequence consisting of closed connected topological
$n$-manifolds $X_i$
and maps $\pi_i:X_{i+1}\to X_i$. Assume furthermore that if the manifolds in
$\cal N$ are oriented then all $X_i$ are also oriented. We say that $\cal J$ is a
{\it weak Jakobsche inverse sequence for $\cal N$} if 
for all $i\ge1$ one can choose finite families ${\cal D}_i$
of collared $n$-disks in $X_i$, partitioned into subfamilies
${\cal D}_i^M:M\in{\cal N}$, such that:

\smallskip
\itemitem{(1)} for each $i\ge1$ the disks in the family 
${\cal D}_i$ are pairwise disjoint;

\itemitem{(2)} for each $i\ge1$ the map $\pi_i$ maps the preimage
$\pi_i^{-1}(X_i\setminus\cup\{ \hbox{int}(D):D\in{\cal D}_i \})$
homeomorphically onto $X_i\setminus\cup\{ \hbox{int}(D):D\in{\cal D}_i \}$;


\itemitem{(3a)} $X_1$ is homeomorphic to one of the manifolds from $\cal N$,
and if the manifolds in $\cal N$ are oriented, we require that this 
homeomorphism respects orientations;

\itemitem{(3b)} for each $i\ge1$, for each $M\in{\cal N}$, and
for any $D\in{\cal D}_i^M$
the preimage $\pi_i^{-1}(D)$ is homeomorphic to $M\setminus\hbox{int}(\Delta)$,
where $\Delta$ is some collared $n$-disk in $M$; furthermore, if the manifolds
in $\cal N$ are oriented, we require that the above homeomorphism
respects the orientations induced from $X_{i+1}$ and from $M$;

\itemitem{(4)} if $i<j$, $D\in{\cal D}_i$, $D'\in{\cal D}_j$, then
$\pi_{i,j}(D')\cap\partial D=\emptyset$, where
$\pi_{i,j}:=\pi_{i}\circ\pi_{i+1}\circ\dots\circ\pi_{j-1}$;

\itemitem{(5)} for each $i\ge1$ the family 
$\{ \pi_{i,j}(D):j\ge i, D\in{\cal D}_j \}$
of subsets of $X_i$ is {\it null}, i.e. the diameters of these subsets
converge to 0; here $\pi_{i,i}$ denotes the identity map
on $X_i$;

\itemitem{(6)} for any $i\ge1$ and each $M\in{\cal N}$ the set
$\bigcup_{j=i}^\infty\pi_{i,j}(\cup{\cal D}_j^M)$
is dense in $X_i$.


\medskip\noindent
{\bf Remarks.}

\item{(1)} It follows from conditions (1), (2), (3a) and (3b)
that each $X_i$ is the connected sum of a family of manifolds
each homeomorphic to one of the manifolds in $\cal N$;
moreover, if the manifolds in $\cal N$ are oriented, the above
mentioned homeomorphisms and the connected sum respect
the orientations.

\item{(2)} In the case when the manifolds in $\cal N$ are oriented,
conditions (1)-(5) in Definition 3.E.3 coincide with conditions
(1)-(6) in [J2], Section 2, p. 82.

\item{(3)} Condition (6) in Definition 3.E.3 implies condition (7)
in [J2], but it is essentially weaker than the conjunction of conditions 
(7) and (8) of [J2] (except when the family $\cal M$ consists of
a single manifold $M$, in which case (6) is equivalent to
the conjunction of (7) and (8), as it was observed and
exploited in [Fi] and [Z1]).

\medskip\noindent
{\bf 3.E.4 Corollary.}
{\it Given ${\cal N}=\{M_1,\dots,M_k  \}$ as in Definition 3.E.3, 
the limit $\lim_{\longleftarrow}{\cal J}$
of any weak 
Jakobsche inverse sequence $\cal J$ for $\cal N$  is homeomorphic to the Jakobsche space
${\cal X}(M_1\#\dots\# M_k)$.}

\medskip\noindent
{\bf Proof:}
First, observe that by conditions (1)-(5) of Definition 3.E.3, 
there is a tree system $\cal M$ of manifolds from $\cal N$ such that
$\cal J$ has the form of an inverse sequence associated to $\cal M$,
as in Example 2.C.7, for an appropriate choice of a conical
family of extended spaces and maps. More precisely, the constituent
spaces of $\cal M$ coincide with the spaces $Y$ of the following two kinds:
\item{(1)} for any $i\ge1$ put 
$$
{\cal D}_i'=\{ D\in{\cal D}_i \hbox{ such that there is
no $j<i$ with $\pi_{i,j}(D)\subset D'$ for some $D'\in{\cal D}_j$} \},
$$
and set 
$Y=X_1\setminus\bigcup_{i=1}^\infty\bigcup_{D\in{\cal D}'_i}\pi_{1,i}(\hbox{int}(D))$;

\item{(2)} for any $m\ge1$, any $\Delta\in{\cal D}_m$, 
and any $i\ge m+1$ put
${\cal D}_{\Delta,i}$ to be the family of all
$D\in{\cal D}_i$ such that $\pi_{m,i}(D)\subset\Delta$ 
and there is
no $m+1\le j<i$ with $\pi_{i,j}(D)\subset D'$ for some $D'\in{\cal D}_j$;
set $Y=\pi_m^{-1}(\Delta)\setminus\bigcup_{i=m+1}^\infty\bigcup_{D\in{\cal D}_{\Delta,i}}\pi_{m+1,i}(\hbox{int}(D))$.

\noindent
We skip further explanations and justifications concerning this first observation, and we note that, due to Theorem 2.B.4, we have
$\lim_{\longleftarrow}{\cal J}=\lim{\cal M}$. 

Next, it follows fairly directly from condition (6) of Definition 3.E.3
that the tree system $\cal M$ of manifolds from $\cal N$, as above, 
is dense 
and weakly saturated. The assertion follows then directly from
Theorem 3.E.2.

\medskip

We finish with presenting briefly a more restrictive, but slightly less technical 
than weak Jakobsche inverse sequence, concept of a {\it special Jakobsche
inverse sequence}. 

\medskip\noindent
{\bf 3.E.5 Definiton.}
Let $\cal N$ be a finite family of closed connected $n$-dimensional 
topological manifolds,
either all oriented, or at least one of which is non-orientable.
Let $${\cal G}=(\{ X_i:i\ge1 \},\{ \pi_i:i\ge1 \})$$
be an inverse sequence consisting of closed connected topological
$n$-manifolds $X_i$
and maps $\pi_i:X_{i+1}\to X_i$. Assume furthermore that if the manifolds in
$\cal N$ are oriented then all $X_i$ are also oriented. We say that $\cal G$ is a
{\it special Jakobsche inverse sequence for $\cal N$} if 
for all $i\ge1$ one can choose finite sets ${\cal Q}_i$
in $X_i$, partitioned into subsets
${\cal Q}_i^M:M\in{\cal N}$, such that:

\smallskip

\itemitem{(1)} for each $i\ge1$ the map $\pi_i$ maps the preimage
$\pi_i^{-1}(X_i\setminus{\cal Q}_i)$
homeomorphically onto $X_i\setminus{\cal Q}_i$;

\itemitem{(2a)} $X_1$ is homeomorphic to one of the manifolds from $\cal N$,
and if the manifolds in $\cal N$ are oriented, we require that this 
homeomorphism respects orientations;

\itemitem{(2b)} for each $i\ge1$, for each $M\in{\cal N}$, and
for any $q\in{\cal Q}_i^M$
the preimage $\pi_i^{-1}(q)$ is a submanifold of $X_{i+1}$ 
homeomorphic to $M\setminus\hbox{int}(\Delta)$,
where $\Delta$ is some collared $n$-disk in $M$; furthermore, if the manifolds
in $\cal N$ are oriented, we require that the above homeomorphism
respects the orientations induced from $X_{i+1}$ and from $M$;

\itemitem{(3)} if $i<j$ and $q\in{\cal Q}_i$, then
${\cal Q}_j\cap\partial[\pi_{ij}^{-1}(q)]=\emptyset$;

\itemitem{(4)} for any $i\ge1$ and each $M\in{\cal N}$ the set
$\bigcup_{j=i}^\infty\pi_{i,j}({\cal Q}_j^M)$
is dense in $X_i$.

\medskip\noindent
{\bf Remark.} 
Condition (3) in the above definition requires a comment. Note that,
by condition (2b), if $q\in{\cal Q}_i$ and if $j=i+1$ then $\pi_{ij}^{-1}(q)$
is an $n$-submanifold with boundary in $X_j$.
Moreover, if for some $j>i$ we have that $\pi_{ij}^{-1}(q)$
is an $n$-submanifold with boundary in $X_j$, and if 
${\cal Q}_j\cap\partial[\pi_{ij}^{-1}(q)]=\emptyset$,
then it follows from conditions (1) and (2b) that
$\pi_{i,j+1}^{-1}(q)$ is an $n$-submanifold with boundary in $X_{j+1}$.
Thus, by induction, all the preimages $\pi_{ij}^{-1}(q)$ occuring
in the statement of condition (3) are $n$-submanifolds with boundary
in the corresponding manifolds $X_j$, and hence it makes sense
to speak of their boundaries $\partial[\pi_{ij}^{-1}(q)]$.

\medskip
An argument similar to that in the proof of Corollary 3.E.4
shows that a special Jakobsche inverse sequence $\cal G$ (for $\cal N$) 
has the form of the standard inverse sequence
associated to a weakly saturated tree system of manifolds
from $\cal N$.
In this argument
the reference  to Theorem 2.B.4 has to be replaced with the corresponding
reference to Proposition 1.D.1. As a consequence, we get the following.

\medskip\noindent
{\bf 3.E.6 Corollary.}
{\it Given ${\cal N}=\{M_1,\dots,M_k  \}$ as in Definition 3.E.5, 
the limit $\lim_{\longleftarrow}{\cal G}$
of any special 
Jakobsche inverse sequence $\cal G$ for $\cal N$  is homeomorphic to the Jakobsche space
${\cal X}(M_1\#\dots\# M_k)$.}

\medskip\noindent
{\bf Remark.}
It is not hard to see that any special Jakobsche inverse sequence
for $\cal N$ is also a weak Jakobsche inverse sequence for $\cal N$,
but we omit the details. (Obviously, this can be used as another 
justification of Corollary 3.E.6.) The converse is not true.


\bigskip
\bigskip
\bigskip
\centerline{\bf References}

\medskip
\itemitem{[AS]} F. Ancel, L. Siebenmann,
{\it The construction of homogeneous homology manifolds},
Abstracts Amer. Math. Soc. 6 (1985), 92.

\itemitem{[Br]} M. Brown, {\it Some applications of an appropximation
theorem for inverse limits}, Proc. Amer. Math. Soc. 11 (1960), 478--481. 

\itemitem{[Ca]} J. W. Cannon, {\it A positional characterization of the
$(n-1)$-dimensional Sierpi\'nski curve in $S^n$ ($n\ne 4$)},
Fund. Math. 79 (1973), 107--112.

\itemitem{[CD]} R. Charney, M. Davis, {\it Strict hyperbolization},
Topology 34 (1995), 329--350.

\itemitem{[Dav]} R. Daverman, Decompositions of Manifolds, 
Academic Press, 1986.

\itemitem{[DJ]} M. Davis, T. Januszkiewicz, 
{\it Hyperbolization of polyhedra}, J. Differential Geometry 34 (1991),
347--388.

\itemitem{[Dr]} A. Dranishnikov, {\it Cohomological dimension of Markov
compacta}, Top. Appl. 154 (2007), 1341--1358.

\itemitem{[DO]} J. Dymara, D. Osajda, {\it Boundaries of right--angled
hyperbolic buildings}, Fund. Math. 197 (2007), 123--165.

\itemitem{[Eng]} R. Engelking, General Topology, PWN - Polish Scientific
Publishers, Warszawa, 1977.

\itemitem{[Fi]} H. Fischer, {\it Boundaries of right--angled Coxeter 
groups with manifold nerves},
Topology 42 (2003), 423--446.

\itemitem{[Fr]} M. Freedman, {\it The topology of four dimensional
manifolds}, J. Differential Geometry 17 (1982), 357--453.

\itemitem{[J1]} W. Jakobsche, {\it The Bing--Borsuk conjecture is stronger 
than the Poincare conjecture},
Fundamenta Mathematicae 106 (1980), 127--134.

\itemitem{[J2]} W. Jakobsche, {\it Homogeneous cohomology manifolds which are inverse limits},
Fundamenta Mathematicae 137 (1991), 81--95.

\itemitem{[JS]} T. Januszkiewicz, J. \'Swi\c{a}tkowski,
{\it Simplicial nonpositive curvature},
Publ. Math. IHES 104 (1) (2006), 1--85.

\itemitem{[KB]}
I. Kapovich, N. Benakli, {\it Boundaries of hyperbolic groups,} 
in: Combinatorial and geometric group theory 
(New York, 2000/Hoboken, NJ, 2001), 39--93,
Contemp. Math. 296, Amer. Math. Soc., Providence, RI, 2002.

\itemitem{[KK]}
M. Kapovich, B. Kleiner, 
{\it Hyperbolic groups with low-dimensional boundary}, 
Ann. Sci. ENS 33 (2000) 647--669.

\itemitem{[P\'S]} P. Przytycki, J. \'Swi\c atkowski,  
{\it Flag-no-square triangulations and Gromov boundaries
in dimension 3}, Groups, Geometry \& Dynamics  3 (2009), 453--468.

\itemitem{[Qu]} F. Quinn, {\it Ends of maps. III: dimensions 4 and 5},
J. Diff. Geom. 17 (1982), 503--521.

\itemitem{[St]} P.R. Stallings, {\it An extension of Jakobsches construction
of $n$--ho\-mo\-geneous continua to the
nonorientable case}, in Continua (with the Houston Problem Book),
ed. H. Cook, W.T. Ingram, K. Kuperberg, A. Lelek, P. Minc,
Lect. Notes in Pure and Appl. Math. vol. 170 (1995), 347--361.

\itemitem{[Sw1]} J. \'Swi\c atkowski, {\it Fundamental pro-groups and Gromov boundaries
of 7-systolic groups,} Journal of the London Mathematical Society 80 (2009), 649--664.

\itemitem{[Sw2]} J. \'Swi\c atkowski, {\it Trees of manifolds as boundaries
of spaces and groups}, preprint, 2013.

\itemitem{[Z1]} P. Zawi\'slak, {\it Trees of manifolds and boundaries of 
systolic groups}, Fund. Math. 207 (2010), 71--99.

\itemitem{[Z2]} P. Zawi\'slak, {\it Trees of manifolds with boundaries},
Colloquium Mathematicum, to appear.

\bigskip

\noindent  Instytut Matematyczny, Uniwersytet Wroc\l awski, 

\noindent pl. Grunwaldzki 2/4, 50-384 Wroc\l aw, Poland

\smallskip
\noindent E-mail: {\tt Jacek.Swiatkowski@math.uni.wroc.pl}

\bye